\documentclass[reqno,11pt,letterpaper]{amsart}
\usepackage[T1]{fontenc}
\usepackage[utf8]{inputenc}
\usepackage{lmodern,amsmath,amssymb,mathtools,mathrsfs}
\usepackage[margin=25mm]{geometry}
\usepackage{microtype,xcolor,enumitem,booktabs,array}
\IfFileExists{dsfont.sty}{\usepackage{dsfont}}{}
\providecommand{\mathds}[1]{\mathbf{#1}}
\definecolor{red}{RGB}{178,28,38}
\colorlet{RED}{red}
\usepackage[hidelinks]{hyperref}
\hypersetup{pdftitle={Dynamics of 3D vortex rings},pdfauthor={Ping Zhang and Yibin Zhang}}

\newcommand{\newcom}{\newcommand}
\newcom{\corr}[1]{{\color{red}#1}}
\newcom{\dd}{\,d}
\newcom{\andf}{\quad\hbox{and}\quad}
\newcom{\with}{\quad\hbox{with}\quad}
\newcom{\al}{\alpha}
\newcom{\eps}{\epsilon}
\newcom{\ve}{\varepsilon}
\newcom{\ga}{\gamma}
\newcom{\Ga}{\Gamma}
\newcom{\ka}{\kappa}
\newcom{\Lam}{\Lambda}
\newcom{\lam}{\lambda}
\newcom{\Sig}{\Sigma}
\newcom{\sig}{\sigma}
\newcom{\tht}{\theta}
\newcom{\om}{\omega}
\newcom{\oo}{\infty}
\newcom{\grad}{\nabla}
\newcom{\lap}{\Delta}

\newcom{\cA}{{\mathcal A}}
\newcom{\cB}{{\mathcal B}}
\newcom{\cC}{{\mathcal C}}
\newcom{\cD}{{\mathcal D}}
\newcom{\cE}{{\mathcal E}}
\newcom{\cF}{{\mathcal F}}
\newcom{\cG}{{\mathcal G}}
\newcom{\cH}{{\mathcal H}}
\newcom{\cI}{{\mathcal I}}
\newcom{\cJ}{{\mathcal J}}
\newcom{\cK}{{\mathcal K}}
\newcom{\cL}{{\mathcal L}}
\newcom{\cM}{{\mathcal M}}
\newcom{\cN}{{\mathcal N}}
\newcom{\cO}{{\mathcal O}}
\newcom{\cP}{{\mathcal P}}
\newcom{\cQ}{{\mathcal Q}}
\newcom{\cR}{{\mathcal R}}
\newcom{\cS}{{\mathcal S}}
\newcom{\cT}{{\mathcal T}}
\newcom{\cU}{{\mathcal U}}
\newcom{\cV}{{\mathcal V}}
\newcom{\cW}{{\mathcal W}}
\newcom{\cX}{{\mathcal X}}
\newcom{\cY}{{\mathcal Y}}
\newcom{\cZ}{{\mathcal Z}}
\newcom{\R}{\mathbb{R}}
\newcom{\T}{\mathbb{T}}
\newcom{\N}{\mathbb{N}}
\newcom{\Z}{\mathbb{Z}}
\newcom{\C}{\mathbb{C}}
\newcom{\E}{\mathbb{E}}
\newcom{\PP}{\mathbb{P}}
\newcom{\HH}{\mathbb{H}}

\newcom{\f}{\frac}
\newcom{\p}{\partial}
\newcom{\pt}{\p_t}
\newcom{\ov}{\overline}
\newcom{\wt}{\widetilde}
\let\wh=\widehat
\def\eqdefa{\buildrel\hbox{\footnotesize def}\over =}
\newcom{\M}{{\scriptscriptstyle (M)}}
\newcom{\MP}{{\scriptscriptstyle (M+1)}}
\newcom{\EL}{{\scriptscriptstyle E}}
\newcom{\NS}{{\scriptscriptstyle NS}}
\newcom{\RR}{{\rm R}}
\newcom{\Ein}{\text{Ein}}
\newcom{\BS}{\mathrm{BS}}
\DeclareMathOperator{\supp}{supp}

\newcom{\br}{{\bar r}}
\newcom{\bz}{{\bar z}}
\newcom{\dbr}{\dot{\br}}
\newcom{\dbz}{\dot{\bz}}
\newcom{\bxi}{\langle\xi\rangle}
\newcom{\ive}{{\mathrm{I}_{\ve}}}
\newcom{\iive}{{\mathrm{II}_{\ve}}}
\newcom{\iiive}{{\mathrm{III}_{\ve}}}

\newtheorem{Theorem}{Theorem}[section]
\newtheorem{Proposition}[Theorem]{Proposition}
\newtheorem{Lemma}[Theorem]{Lemma}

\theoremstyle{remark}
\newtheorem{Remark}[Theorem]{Remark}

\renewcommand{\theequation}{\thesection.\arabic{equation}}
\allowdisplaybreaks[2]
\setlist[enumerate]{leftmargin=*,itemsep=3pt,topsep=4pt}
\setlist[itemize]{leftmargin=*,itemsep=3pt,topsep=4pt}
\title[Dynamics of a 3D vortex ring]{Dynamics of a single 3D vortex ring}
\author[P. Zhang]{Ping Zhang}
\address{State Key Laboratory of Mathematical Sciences, Academy of Mathematics and Systems Science, Chinese Academy of Sciences, Beijing 100190, China; and School of Mathematical Sciences, University of Chinese Academy of Sciences, Beijing 100049, China}
\email{zp@amss.ac.cn}
\author[Y. Zhang]{Yibin Zhang}
\address{Academy of Mathematics and Systems Science, Chinese Academy of Sciences, Beijing 100190, China}
\email{zhangyibin22@mails.ucas.ac.cn}
\date\today
\begin{document}

\begin{abstract}
We study the axisymmetric Navier--Stokes flow without swirl originating from a circular vortex ring. For each fixed $s\in(0,1)$ and sufficiently large circulation Reynolds number
$\Gamma/\nu$, we prove that, on the time scale $t\lesssim T_{\rm adv}(\Gamma/\nu)^{s}$, the
vorticity $\omega_\theta$ is approximated by a traveling wave solution featuring a Gaussian-type core with scale $\cO(\sqrt{\nu t})$, whose off-axis distance remains constant and which continuously translates along the axial direction. In particular, the leading-order axial translation velocity is given by the Kelvin-Saffman formula. The proof relies on the approximation and Arnold energy framework of Gallay and \v{S}ver\'ak \cite{GS2024}. We modify the weight used in \cite{GS2024} and estimate the relative-vorticity error to obtain favorable boundary dissipation at the symmetry axis and close the energy estimate.
\end{abstract}
\maketitle


\color{black}
\renewcommand{\leftmark}{P. ZHANG AND Y. ZHANG}
\renewcommand{\rightmark}{DYNAMICS OF A SINGLE 3D VORTEX RING}

\renewcommand{\theequation}{\thesection.\arabic{equation}}
\setcounter{equation}{0}

\section{Introduction}

Coherent vortex structures play a central role in incompressible fluid
dynamics. A vortex ring is a localized toroidal distribution of vorticity whose induced velocity transports it along its axis, and  it provides a basic model of such coherent structures~\cite{SL1992}. Vortex rings arise in starting jets and biological propulsion, where vortex formation is closely related to momentum transfer~\cite{KG2023,D2009}. These applications motivate understanding how a concentrated ring
retains its structure during propagation. At high Reynolds numbers, the core can remain thin while the ring
travels over distances much larger than its radius. A basic problem
is therefore to determine how long an asymptotic core profile describes
the flow and how core deformation affects its trajectory. Circular
rings provide a natural setting for this question: axisymmetry
preserves their circular geometry, while curvature generates a
self-induced translation that must be tracked together with viscous
spreading. We study the free evolution of a single ring using the
three-dimensional Navier--Stokes equations in the axisymmetric class
without swirl:
\begin{equation}\label{eq 1.1}
\begin{cases}
\displaystyle
\partial_t\omega_\theta+\left\{\varphi,\frac{\omega_\theta}{r}\right\}
=\nu\left((\partial_r^2+\partial_z^2)\omega_\theta
+\partial_r\left(\frac{\omega_\theta}{r}\right)\right),\\
\omega_\theta|_{r=0}=0, \quad (r,z)\in\HH \eqdefa (0,\infty)\times\R.
\end{cases}
\end{equation}
Here the Poisson bracket $\{f,g\}\eqdefa \partial_r f\,\partial_z g-\partial_z f\,\partial_r g$,  $\omega_\theta$ denotes the azimuthal component of vorticity, $\nu>0$ is the kinematic viscosity, and the absence of swirl means that the azimuthal velocity is zero. The stream function $\varphi$  satisfies the following elliptic equation:
\begin{equation}\label{eq 1.2}
\begin{cases}
\displaystyle
    -\partial_r(r^{-1}\partial_r\varphi) - \partial_z(r^{-1}\partial_z \varphi)=\omega_\theta,  \\
     \varphi|_{r=0}=\partial_r\varphi|_{r=0}=0, \quad (r,z)\in \HH,
\end{cases}
\end{equation}
 with the associated
velocity vanishing at spatial infinity
\begin{equation}\label{eq 1.3}
 u_r=-\frac1r\partial_z\varphi \rightarrow 0,
\quad u_z=\frac1r\partial_r\varphi \rightarrow 0, \qquad \text{as}\, (r,z)\rightarrow \oo.
\end{equation}
We consider the solution to \eqref{eq 1.1} with initial data
\begin{equation}\label{eq 1.4}
\omega_\theta(0)=\Gamma\delta_{(r_0,z_0)}.
\end{equation}
Here $r_0>0$ is the initial radius, $z_0\in\R$ denotes the initial axial position, and $\Gamma>0$ is the circulation parameter. The opposite sign can be treated by reversing the axial orientation. At positive times, viscosity regularizes this mass into a core of characteristic width $\sqrt{\nu t}$. Locally, the leading profile is Gaussian, as for a straight viscous filament. The circular geometry nevertheless generates an axial velocity whose leading term depends logarithmically on the ratio between the ring radius and the core width. Thus even when the core remains thin, the ring can travel over distances much larger than its radius. A long-time approximation must capture both this translation and the evolving internal structure.

To be more specific, we denote the advection time, the diffusion time, and the inverse circulation Reynolds number by
\begin{equation}\label{eq 1.5}
T_{\mathrm{adv}}=\frac{r_0^2}{\Gamma},\qquad
T_{\mathrm{dif}}=\frac{r_0^2}{\nu},\qquad
\delta=\frac{\nu}{\Gamma},\qquad
T_{\mathrm{dif}}=\delta^{-1}T_{\mathrm{adv}}.
\end{equation}
In the moving coordinates introduced below, the leading axial velocity is given by the Kelvin-Saffman formula:
\begin{equation}\label{eq 1.6}
V(t)=\frac{\Gamma}{4\pi r_0}
\left(\log\frac{ r_0}{\sqrt{\nu t}}
+\frac32\log 2+\frac{\gamma_E-1}{2}\right),
\end{equation}
where $\gamma_E$ denotes Euler's constant. Set $\widetilde t=\Gamma t$, $\widetilde\omega=\omega_\theta/\Gamma$,
$\widetilde\varphi=\varphi/\Gamma$, and $\widetilde\nu=\nu/\Gamma$, and then drop the tildes, we assume throughout the rest of the paper that $\Gamma=1$. Thus $\nu$ is the inverse circulation Reynolds number, and we study the regime $\nu\ll1$.

Our aim is to justify this description on time scales much longer than the advection time and approaching, in algebraic order, the diffusion time. For the reader's convenience, we first recall the relevant well-posedness and approximation results, and then introduce the formulation in which our main theorem is stated.

\subsection{Classical vortex rings and viscous asymptotics}

The mathematical study of vortex rings includes both the construction of steadily translating inviscid rings and the evolution of viscous cores. In the Euler setting, Fraenkel \cite{Fraenkel1970} established the existence of steady rings of small cross-section and obtained asymptotic descriptions of their propagation speed and shape. Norbury \cite{Norbury1973} gave a numerical description of a family of steady axisymmetric rings connecting the thin-ring regime with Hill's spherical vortex. Friedman and Turkington \cite{FT1981} developed a variational construction of steady vortex rings and derived thin-core asymptotic estimates. Together, these works provide a mathematical account of the relation between core geometry and self-induced motion for inviscid traveling profiles.

Viscosity introduces a further question, because the core spreads and the propagation speed varies with time. Saffman \cite{Saffman1970} extended the calculation of ring velocity to an unsteady viscous core. For a Gaussian core, the resulting Kelvin--Saffman formula has the logarithmic dependence displayed in \eqref{eq 1.6}. Fukumoto and Moffatt \cite{FM2000} subsequently developed a higher-order matched asymptotic expansion, accounting for the elliptical deformation induced by the ring's own strain and for radial evolution. Their analysis shows how corrections to the core enter the translation speed and the geometry of a thin viscous ring.

These steady-ring constructions and viscous asymptotic calculations identify the profiles and motion laws that a rigorous dynamical approximation should recover. The remaining issue is to control the difference between a solution of the singular Navier--Stokes initial-value problem and a suitably corrected profile over a time interval that grows as the viscosity decreases. This requires both an accurate core expansion and stability estimates that remain effective in the large circulation Reynolds number regime.

\subsection{Singular initial data and rigorous long-time dynamics}

The singular initial-value problem for circular vortex filaments is partially understood in the axisymmetric class without swirl. Feng and \v{S}ver\'ak \cite{FS2015} constructed global solutions originating from circular filaments of arbitrary strength under a sign condition. In scale-invariant vorticity spaces, Gallay and \v{S}ver\'ak \cite{GS2015} proved global well-posedness for integrable vorticity and treated finite-measure data whose atomic part is sufficiently small relative to the viscosity. For a single circular filament with arbitrary circulation Reynolds number, Gallay and \v{S}ver\'ak \cite{GS2019} established uniqueness and the leading Gaussian core asymptotics at small positive times. L\'evy and Liu \cite{LL2018} extended uniqueness to positive combinations of circular vortex filaments. Vortex filament solutions associated with more general smooth curves were constructed by Bedrossian, Germain, and Harrop-Griffiths \cite{BGH2023}, in a setting where the geometry is no longer fixed by axisymmetry.

A related line of research studies initially smooth vorticity concentrated in thin rings. Butt\`a, Cavallaro, and Marchioro \cite{BCM2022} justified concentration and axial translation for several separated viscous rings in a joint limit of vanishing initial thickness and viscosity, with appropriately scaled circulations. Their result allows several rings and general concentrated initial cores. The initial thickness, the viscosity, and the circulation scaling are linked in that limit. Here the initial core is a Dirac mass and the circulation is fixed, so that its positive-time width is generated by viscosity itself.

For this singular-filament problem, Gallay and \v{S}ver\'ak \cite{GS2024} constructed a refined vortex-core approximation and controlled its error with a quadratic energy inspired by Arnold's variational method. Their theorem gives a time interval of order $T_{\mathrm{adv}}\delta^{-\sigma}$ for a sufficiently small exponent $\sigma\in(0,1/3)$. Their more precise approximation incorporates corrections to the Kelvin--Saffman formula and gives stronger error control than the leading translated heat profile. This work provides the starting point for our analysis. Related long-time approximation methods for planar vortex pairs were developed by Dolce and Gallay \cite{DG2026} and Zhang and Zhang \cite{ZZ26}, combining a moving frame, higher-order cores, and an energy adapted to the Euler transport structure.

A small error in the translation speed can accumulate over a long interval, and the core approximation must remain stable under a transport operator that becomes large as the viscosity tends to zero. The central question is therefore whether the single-ring description can be maintained on every fixed algebraic time scale below diffusion. The theorem proved here answers this question in the axisymmetric setting by establishing a Gaussian approximation, with a corrected trajectory.

For a vortex ring, extending the approximation toward the diffusion time requires keeping track of several features of the axisymmetric geometry. The rescaled domain has a moving boundary corresponding to the symmetry axis; the  vorticity loses mass through that axis; and the Biot--Savart operator changes with the aspect ratio of the core. In addition, logarithmic coefficients and the time dependence of the radius enter the higher-order construction. Our analysis incorporates these effects into the approximate dynamics, the energy identity, and the modulation estimates. This leads to a Gaussian approximation valid for every fixed algebraic time exponent below the diffusion scale.

\subsection{Main result}

Our main theorem describes the vorticity and the axial motion of the ring throughout every fixed algebraic time interval below the diffusion scale. Recall that we assume $\Ga=1$ throughout this paper.

\begin{Theorem}\label{s1:main}
Let $s\in(0,1)$ and $(r_0,z_0)\in\HH$. There exist $\nu_0=\nu_0(r_0,s)>0$, $C=C(r_0,s)>0$, and an integer $N\ge0$ such that, for every $0<\nu\le\nu_0$, the unique axisymmetric mild solution of \eqref{eq 1.1} with initial trace $\omega_\theta(0)=\delta_{(r_0,z_0)}$ satisfies
\begin{equation}\label{eq 1.7}
\begin{aligned}
\int_{\HH}\left|\omega_\theta(t,r,z)
-\frac1{4\pi\nu t}\exp\left(-\frac{(r-r_0)^2+(z-\bar z_a(t))^2}{4\nu t}\right)\right|
\dd r\dd z
\le C\sqrt{\nu t}\,|\log(\nu t)|^N,
\end{aligned}
\end{equation}
uniformly for $0<t\le\nu^{-s}$. The approximate trajectory satisfies $\bar z_a(0)=z_0$ and
\[
\dot{\bar z}_a(t)=\frac1{4\pi r_0}
\left(\log\frac{r_0}{\sqrt{\nu t}}
+\frac32\log2+\frac{\gamma_E-1}{2}\right)
+O\bigl((\nu t)^{\frac12-}\bigr).
\]
\end{Theorem}

\begin{Remark}
  Notice that the leading-order term of $\dbz_a(t)$ is given by the Kelvin--Saffman formula \eqref{eq 1.6}.
 The notation $\frac12-$ means $\frac12-\eta$ for an arbitrarily small fixed $\eta>0$, with the implicit constant allowed to depend on $\eta$. 
\end{Remark}

\begin{Remark}
    The bound $\cO(\sqrt{\nu t} |\log(\nu t)|^N)$ can be sharpened by incorporating higher-order corrections. And the solution $\omega_\theta$ exhibits more refined localization properties than the $L^1$ bound. Roughly speaking, $\omega_\theta$ remains highly concentrated near the trajectory $(r_0,\bz_a(t))$, while away from this point it decays rapidly at the rate $\exp(-\Big(\frac{(r-r_0)^2+(z-\bz_a(t))^2}{4\nu t}\Big)^\gamma)$ for some $\ga\ll1$, in $L^2$ integral sense.
\end{Remark}

\subsection{Proof strategy and new ingredients}

The proof builds on the approximation and Arnold-energy framework of \cite{GS2024}, together with the long-time constructions in \cite{DG2026,ZZ26}. We organize the argument into three steps.

\smallskip
\noindent\textbf{Step 1. Higher-order approximate cores and the functional relation.}
Section~\ref{s2:section} collects the expansion of the Biot--Savart operator and the properties of the planar linearized operators. In Section~\ref{s3:section}, for each fixed integer $M\ge2$, we construct an approximate profile $\Omega_a^{\M}$ and frame velocities $(\dot{\bar r}_a^{\M}$, $\dot{\bar z}_a^{\M})$ so that the residual satisfies
\[
\RR\bigl[\Omega_a^{\M},\dot{\bar r}_a^{\M},\dot{\bar z}_a^{\M}\bigr]
=O_{\cZ}^{L}\bigl(\nu^{-1}\mu^{M+1}+\nu\mu\bigr).
\]
However, in our construction, due to the logarithmic corrections, the time dependence of the coefficients is retained in the induction. Thus at order $\mu^n$, the coefficient equation contains the operator $t\partial_t+n/2-\cL$, the inversion properties of which are established in Lemma \ref{s2:lemma-inversion}.

Propositions~\ref{prop 3.2} and~\ref{prop 3.3} then establish, to the required order, a local relation between
\[
\Psi_a\eqdefa \Phi_a-\frac{\bar r_a\dot{\bar z}_a}{2}(\varpi^2-1),\qquad Q_a\eqdefa \frac{\Omega_a}{\varpi}.
\]
This relation supplies the core weight in Arnold's energy and makes the large linearized transport term compatible with the stability estimate.

\smallskip
\noindent\textbf{Step 2. An adapted energy and the axisymmetric cancellations.} This step differs from the strategy in \cite{GS2024}. We measure the error in relative vorticity,
\[
q\eqdefa \varpi^{-1}\Omega-Q_a,\qquad
\phi\eqdefa \BS_\varepsilon[\varpi q],
\]
and derive its exact evolution equation in Section~\ref{s4:section}. The corresponding energy and dissipation functionals are
\begin{equation*}
    \begin{split}
    &E_\varepsilon[q]
=\frac12\int_{\HH_\varepsilon}W_\varepsilon q^2\dd\xi
-\frac12\int_{\HH_\varepsilon}\phi q\varpi\dd\xi,
\qquad W_\varepsilon=\varpi H_\varepsilon,\\
&\|q\|_{D_\ve}^2\eqdefa
       \int_{\HH_\ve}W_\ve
        (|\grad q|^2+\rho_\ga^2q^2+q^2)\dd\xi
       +\ve\int_{\p\HH_\ve}H_\ve q^2\dd\xi_2.
    \end{split}
\end{equation*}
Here in terms of the vorticity error $\eta=\varpi q$, its quadratic weight is $H_\ve/\varpi$, which is singular at the symmetry axis, in contrast to the weight in \cite{GS2024}. Controlling this weight requires weighted Hardy inequalities and stream-function estimates that retain the  axis factors.

In the core, $H_\varepsilon$ is obtained by differentiating the approximate functional relation. It is constant in the intermediate region and grows at a stretched-exponential rate in the far field. This choice connects the Euler structure near the core with the localization needed away from it. Two cancellations are essential to the energy estimate. Since $H_\varepsilon$ is constant in the intermediate region, the contribution involving its spatial derivative vanishes there. Moreover, the leading time derivative of the Biot--Savart energy cancels a term generated by diffusion. The remaining boundary trace term at the symmetry axis contributes with a favorable sign to the dissipation. These terms are combined in the exact identity before they are estimated, yielding coercivity and dissipation bounds compatible with the factor $\nu^{-1}$ in the transport operator.

\smallskip
\noindent\textbf{Step 3. Mass, impulse, and axial modulation.}
The energy is coercive after the mass and first-moment components of the error have been controlled. These quantities are estimated using the physical balances in Section~\ref{s5:section}. The  mass has a boundary loss, and a cutoff in the physical radial variable shows that this loss is smaller than any algebraic power of $\ve$. The conserved axial impulse,
\[
\int_{\HH}r^2\omega_\theta(t,r,z)\dd r\dd z=r_0^2,
\]
controls the radial first moment. Then by choosing the free parameter $\bar z(t)$ appropriately, we impose
$M_2[\Omega]=0$. Combining these estimates with the energy inequality closes the bootstrap argument on $0<t\le\nu^{-s}$ for any $s\in(0,1)$ in Section~\ref{s6:section}.

\subsection{Notation}

We use $\beta_a=\log(1/a)$ and $\langle\xi\rangle=(1+|\xi|^2)^{1/2}$. The Gaussian profile, the planar Fokker--Planck operator, the linearized Euler operator, and the Gaussian weighted space are
\[
\begin{gathered}
G(\xi)\eqdefa \frac1{4\pi}e^{-|\xi|^2/4},\qquad
\cL\eqdefa\Delta_\xi+\frac12\xi\cdot\nabla_\xi+1,\\
\Lambda\eqdefa\{\Upsilon,\,\cdot\,\}
+\{\Delta_\xi^{-1}[\,\cdot\,],G\},\qquad
\cY\eqdefa L^2(\R^2;G^{-1}\dd\xi),
\end{gathered}
\]
where
\[
\begin{aligned}
\Upsilon(\xi)=\Delta_\xi^{-1}G
&=\frac1{2\pi}\int_{\R^2}\log|\xi-\eta|\,G(\eta)\dd\eta=\frac1{4\pi}\left(\Ein(|\xi|^2/4)+\log4-\gamma_E\right),
\end{aligned}
\]
with $\Ein(x)=\int_0^x\frac{1-e^{-t}}{t}\dd t.$

The space $\cS_\star(\R^2)$ consists of the smooth multipliers of $\cS(\R^2)$ and $\cS'(\R^2)$: for every multi-index $\alpha\in\mathbb N^2$, there exist $C_\alpha>0$ and $N_\alpha\ge0$ such that
\[
|\partial^\alpha f(\xi)|\le C_\alpha\langle\xi\rangle^{N_\alpha}.
\]
We also define
\[
\cZ=G\cS_\star
=\{f\in\cY:f=Gh\text{ for some }h\in\cS_\star\},
\]
and
\[
\cK=\{F\in C^\infty(\R^+):
F(G)\in\cS_\star(\R^2)\}.
\]
For a function on $\HH_\varepsilon$, membership in $\cZ$ or $\cS_\star$ means that it is a restriction of some functions on $\R^2$ in the corresponding functional space.

For an $\varepsilon$-dependent family, we write $f=O_{\cS_\star}^{L}(\kappa)$ if, for every multi-index $\alpha$, there exist $C>0$ and integers $N,M\ge0$, independent of the small parameters, such that
\[
|\partial^\alpha f(\xi)|
\le C\langle\xi\rangle^N\kappa|\log\varepsilon|^M.
\]
We write $f=O_{\cZ}^{L}(\kappa)$ when $G^{-1}f=O_{\cS_\star}^{L}(\kappa)$ and $F\in\cO_\cK^L(\kappa)$ when $F(G)\in \cO_{\cS_\star}^L(\kappa)$. A function is called $\xi_2$-odd or $\xi_2$-even according as
\[
f(\xi_1,-\xi_2)=-f(\xi_1,\xi_2)
\quad\text{or}\quad
f(\xi_1,-\xi_2)=f(\xi_1,\xi_2).
\]
Finally, $\cP_0$ denotes the projection onto the radial component. For functions on $\R^2$, we define
\begin{equation*}
\bar{m} [f]=\int_{\R^2}f\dd\xi,
\qquad \bar{M}_j[f]=\int_{\R^2}\xi_j f\dd\xi,
\qquad j=1,2.
\end{equation*}

\renewcommand{\theequation}{\thesection.\arabic{equation}}
\setcounter{equation}{0}

\section{Vorticity formulation and basic operators}\label{s2:section}
In this section, we  introduce suitable self-similar coordinates and summarize the properties of   $\BS_\ve$, $\cL$, $\Lam$ and $t\p_t+\ka-\cL$.

\subsection{Vorticity formulation and moving coordinates}
We work with the unique axisymmetric mild solution of \eqref{eq 1.1}, which is smooth for $t>0$ and has the prescribed measure as its initial trace; see \cite{FS2015, GS2019}. 
The stream function selected by the axisymmetric Biot-Savart law is
\begin{equation}\label{eq 2.1}
\varphi(r,z)=\frac1{2\pi}\int_{\HH}\sqrt{rr'}\,
F\left(\frac{(r-r')^2+(z-z')^2}{rr'}\right)
\omega_\theta(r',z')\dd r'\dd z',
\end{equation}
where
\begin{equation}\label{eq 2.2}
\begin{aligned}
F(s)&=\int_0^{\pi/2}
\frac{1-2\sin^2\tau}{\sqrt{\sin^2\tau+s/4}}\dd\tau=
\begin{cases}
\log(8/\sqrt{s})-2+O(s|\log s|),&s\downarrow0,\\[2pt]
\dfrac{\pi}{2s^{3/2}}+O(s^{-5/2}),&s\to\infty.
\end{cases}
\end{aligned}
\end{equation}

To separate the motion of the ring from the spreading of its core, let $(\bar r(t),\bar z(t))\in\HH$ be a differentiable reference point and define
\begin{equation}\label{eq 2.3}
\mu\eqdefa \sqrt{\nu t},\quad \varepsilon\eqdefa \frac{\sqrt{\nu t}}{\bar r(t)},\quad
\xi_1\eqdefa \frac{r-\bar r(t)}{\sqrt{\nu t}},\quad
\xi_2\eqdefa\frac{z-\bar z(t)}{\sqrt{\nu t}},\quad
\varpi(\xi)\eqdefa 1+\varepsilon\xi_1.
\end{equation}
The variable $\varepsilon$ measures the aspect ratio of the core. Since $r=\bar r(t)\varpi(\xi)$, the rescaled domain is
\begin{equation}\label{eq 2.4}
\HH_\varepsilon=\{\xi\in\R^2:\varpi(\xi)>0\}.
\end{equation}
And we define the self-similar profiles of vorticity and stream function via
\begin{equation}\label{eq 2.5}
\begin{aligned}
\Omega(t,\xi)&=\nu t\,\omega_\theta
\bigl(t,\bar r(t)+\sqrt{\nu t}\xi_1,\bar z(t)+\sqrt{\nu t}\xi_2\bigr),\\
\Phi(t,\xi)&=\frac1{\bar r(t)}\varphi
\bigl(t,\bar r(t)+\sqrt{\nu t}\xi_1,\bar z(t)+\sqrt{\nu t}\xi_2\bigr).
\end{aligned}
\end{equation}
They satisfy the elliptic relation
\begin{equation}\label{eq 2.6}
-\nabla_\xi\cdot(\varpi^{-1}\nabla_\xi\Phi)=\Omega,
\end{equation}
with axis conditions
\begin{equation}\label{eq 2.7}
\Omega=0,\qquad \Phi=\partial_1\Phi=0
\qquad\text{on }\partial\HH_\varepsilon.
\end{equation}
Furthermore, by \eqref{eq 2.1}-\eqref{eq 2.2}, the elliptic relation can also be written as
\begin{equation}\label{eq 2.8}
\Phi(t,\xi)=\BS_\varepsilon[\Omega](\xi)
=\int_{\HH_\varepsilon}K(\xi,\zeta)\Omega(t,\zeta)\dd\zeta,
\end{equation}
where
\begin{equation}\label{eq 2.9}
K(\xi,\zeta)=\frac{\sqrt{\varpi(\xi)\varpi(\zeta)}}{2\pi}
F\left(\frac{\varepsilon^2|\xi-\zeta|^2}
{\varpi(\xi)\varpi(\zeta)}\right).
\end{equation}

Introducing the Fokker--Planck operator
\[
\cL=\Delta_\xi+\frac12\xi\cdot\nabla_\xi+1,
\]
we obtain the rescaled evolution equation
\begin{equation}\label{eq 2.10}
(t\partial_t-\cL)\Omega
-\varepsilon\partial_1\left(\frac{\Omega}{\varpi}\right)
-\frac{\varepsilon\bar r\dot{\bar r}}{\nu}\partial_1\Omega
+\nu^{-1}\left\{\Phi-\frac{\bar r\dot{\bar z}}2\varpi^2,
\frac{\Omega}{\varpi}\right\}_\xi=0.
\end{equation}
Here $\{f,g\}_\xi\eqdefa \p_{\xi_1}f \p_{\xi_2}g-\p_{\xi_1}g \p_{\xi_2}f $. In what follows, Poisson brackets between rescaled quantities act on $\xi$, and we omit the subscript $\xi$. In this formulation, the main analytical difficulty is visible in the factor $\nu^{-1}$ multiplying the transport term in \eqref{eq 2.10}. Uniform estimates therefore require an energy that reflects the structure of this transport operator, together with control of the moving boundary and the time-dependent kernel.

\subsection{Mass and first moments}

We record the identities that determine the frame velocities and enter the stability argument. For an integrable profile with finite first moments, define
\begin{equation}\label{eq 2.11}
m[\Omega]=\int_{\HH_\varepsilon}\Omega\dd\xi,
\qquad M_j[\Omega]=\int_{\HH_\varepsilon}\xi_j\Omega\dd\xi,
\qquad j=1,2.
\end{equation}

\begin{Lemma}\label{lem 1.4}
Let $\Omega$ satisfy \eqref{eq 2.6}--\eqref{eq 2.10}. Then
\begin{equation}\label{eq 2.12}
t\frac{\mathrm{d}}{\mathrm{d}t}m[\Omega]
=-2\int_\R\partial_1\Omega(t,-\varepsilon^{-1},\xi_2)\dd\xi_2,
\end{equation}
\begin{equation}\label{eq 2.13}
\begin{aligned}
\left(t\frac{\mathrm{d}}{\mathrm{d}t}+\frac12\right)M_1[\Omega]
={}&-\nu^{-1}\int_{\HH_\varepsilon}\varpi^{-1}\partial_2\Phi\,\Omega\dd\xi
-\varepsilon\int_{\HH_\varepsilon}\varpi^{-1}\Omega\dd\xi\\
&+\frac2\varepsilon\int_\R\partial_1\Omega(t,-\varepsilon^{-1},\xi_2)\dd\xi_2
-\frac{\varepsilon\bar r\dot{\bar r}}{\nu}m[\Omega],
\end{aligned}
\end{equation}
and
\begin{equation}\label{eq 2.14}
\begin{aligned}
\left(t\frac{\mathrm{d}}{\mathrm{d}t}+\frac12\right)M_2[\Omega]
={}&\nu^{-1}\int_{\HH_\varepsilon}\varpi^{-1}\partial_1\Phi\,\Omega\dd\xi
-\frac{\varepsilon\bar r\dot{\bar z}}{\nu}m[\Omega]\\
&-2\int_\R\xi_2\partial_1\Omega(t,-\varepsilon^{-1},\xi_2)\dd\xi_2.
\end{aligned}
\end{equation}
\end{Lemma}
\begin{proof}
    See Appendix~\ref{appA}.
\end{proof}

\subsection{Expansion of the Biot--Savart law}
The scalar kernel $F$ has the convergent expansion
\[
 F(s)=\log\f8{\sqrt s}\sum_{n=0}^{\oo}c_ns^n
                      +\sum_{n=0}^{\oo}d_ns^n,\qquad 0<s<4,
\]
where
\[
 c_0=1,\qquad d_0=-2,\qquad c_1=\f3{16},\qquad d_1=-\f1{16}.
\]
Let $D=|\xi-\zeta|$ and $L=\log(8/D)$. The kernel of
$\BS_\ve$ consequently has the following formal expansion (see \cite[Lemma~3.3]{GS2024}),
\begin{equation}\label{eq 2.15}
 K(\xi,\zeta)
  =(\beta_\ve+L)\sum_{m=0}^{\oo}\ve^mP_m
                         +\sum_{m=0}^{\oo}\ve^mQ_m,
\end{equation}
where $P_m,Q_m$ are polynomials in $\xi_1,\zeta_1$, and $\xi_2-\zeta_2$.
In particular,
\[
 P_0=\f1{2\pi},\qquad Q_0=-\f1\pi,\qquad
 P_1=\f{\xi_1+\zeta_1}{4\pi},\qquad
 Q_1=-\f{\xi_1+\zeta_1}{4\pi}.
\]
Define
\[
 \BS_m[f](\xi) \eqdefa \beta_\ve P_m[f]+(LP_m)[f] + Q_m[f]  
  \eqdefa \int_{\R^2}\bigl((\beta_\ve+L)P_m+Q_m\bigr)(\xi,\zeta) f(\zeta)\dd\zeta.
\]
At order zero, the expansion gives
\begin{equation}\label{eq 2.16}
 \BS_0[f]
  =-\lap_\xi^{-1}f+
      \f{\beta_\ve+\log8-2}{2\pi}\bar{m}[f].
\end{equation}
Here $\lap_\xi^{-1}f=(2\pi)^{-1}\int_{\R^2}
\log|\xi-\zeta|f(\zeta)\dd\zeta$, and moments of whole-plane coefficients
are integrated over $\R^2$. The following lemmas give basic estimates for the Biot-Savart operator that retain the axis factor $\varpi$. Their proofs are postponed to Appendix~\ref{appA}.

\begin{Lemma}\label{lem 2.2}
Let $0<\ve\ll 1$, $2<p<\oo$ and $\phi=\BS_\ve[\eta]$,  there hold
\begin{align}
 &\|\varpi^{-1}\grad_\xi\phi\|_{L^p}
 \lesssim_p\|\bxi^2\eta\|_{L^2},\label{eq 2.17}\\
 &\|\bxi\varpi^{-1}\grad_\xi\phi\|_{L^\oo}
 \lesssim \|\bxi^2\eta\|_{L^2}^{1/2}
 \bigl(\|\bxi^2\eta\|_{L^2}
       +\|\bxi^2\grad_\xi\eta\|_{L^2}\bigr)^{1/2},\label{eq 2.18}\\
&|\phi(\xi)|\lesssim \sqrt{\varpi(\xi)}
 \bigl(\beta_\ve+\log\bxi\bigr)\|\bxi^2 \eta\|_{L^2}.\label{eq 2.19}
\end{align}
\end{Lemma}
\begin{proof}
   See Appendix \ref{appA}.
\end{proof}

\begin{Lemma}\label{lem 2.1}
Let $0< \ve\ll1$. For any $f=\cO_{\cZ}(1)$ and any integer $N\geq0$, we have 
\begin{equation}\label{eq 2.20}
 \left(\BS_\ve-
                \sum_{m=0}^N\ve^m\BS_m\right)[f]
  =\cO_{\cS_\star}
                 \bigl(\beta_\ve\ve^{N+1}\bigr).
\end{equation}
In particular, if $\supp f \subset \{|\xi|\leq 2\ve^{-\sig_1}\}$ for some $\sig_1 < 1$, then for $|\xi|\leq 2\ve^{-\sig_1}$ and some $p>0$,
\begin{equation}\label{eq 2.21}
    \Big| (\BS_\ve-\BS_0)[f]\Big|+\Big|\varpi^{-1}\grad \BS_\ve[f]-\grad\BS_0[f]\Big| \lesssim \ve\beta_\ve \bxi^p \|\bxi^p f\|_{L^2}.
\end{equation}
\end{Lemma}
\begin{proof}
    See Appendix \ref{appA}.
\end{proof}

\subsection{Properties of $\Lam$ and $\cL$.}
Let $\xi=\rho (\cos\vartheta, \sin\vartheta)$. Expanding $f\in\mathcal Y$ in a Fourier series in the angular variable $\vartheta$, we obtain the
orthogonal decomposition
\begin{equation*}
    \cY= \bigoplus_{n=0}^\oo \cY_n,  \quad \text{with} \ \cY_n \eqdefa \bigl\{ f\in \cY \ :\ f=a(\rho) \cos(n\vartheta) + b(\rho) \sin(n\vartheta)  \bigr\}.
\end{equation*}
We first recall the following basic properties
\begin{Proposition}[\cite{DG2026,ZZ26}]\label{Prop a1}
 {\sl   There hold
\begin{enumerate}
    \item[(1)]  $\cL$ is self-adjoint in $\cY$ with compact resolvent and purely discrete spectrum:
    \begin{equation*}
        \sig(\cL) = \left\{-\frac{n}{2} \ : \ n\in \N\right\}.
    \end{equation*}
    Moreover, $\ker\cL= {\rm span}\{G\}$, $\ker (\f12+\cL) = {\rm span}\{\p_1 G,\p_2G\}$. And the eigenvalue $-\frac{n}{2}$ has multiplicity $n+1$, with eigenfunctions being Hermite functions of degree $n$.


    \item[(2)] $\Lam$ is skew-adjoint in $\cY$ with
    \begin{equation*}
        \ker\Lam = \cY_0 \oplus {\rm span}\bigl\{\p_1 G, \p_2 G\bigr\}.
    \end{equation*}

    \item[(3)] Let
    \begin{equation*}
        \cY'_1 \eqdefa \cY_1 \cap (\ker{\Lam})^\perp = \left\{ f\in \cY_1 \ : \ \int_{\R^2} \xi_j f(\xi)\, d\xi=0 ,\  j=1,2   \right\}.
    \end{equation*}
    Then for $f\in \cY_n \cap \cZ$ if $n\geq 2$,  or $f\in \cY'_1\cap \cZ$ if $n=1$, there exists a unique $\Omega $ in the same space so that  $\Lam \Omega = f$.
     Furthermore, if $f=b(\rho) \sin(n\vartheta)$ (respectively, $f=b(\rho)\cos(n\vartheta)$), then $\Omega= w(\rho)\cos(n\vartheta)$ (respectively, $\Omega=-w(\rho)\sin(n\vartheta)$), where
\begin{equation}\label{eq 2.22}
    w(\rho)\eqdefa - \frac{\rho^2/4}{e^{\rho^2/4}-1} \varphi(\rho) - \frac{2\pi\rho^2 }{n(1-e^{-\rho^2/4})} b(\rho),
\end{equation}
with $\varphi$ solving
\begin{equation}\label{eq 2.23}
    -\varphi'' - \frac{1}{\rho} \varphi' +\left(\frac{n^2}{\rho^2} - \frac{\rho^2/4}{e^{\rho^2/4}-1} \right) \varphi = \frac{2\pi\rho^2 }{n(1-e^{-\rho^2/4})} b(\rho), \quad \rho>0.
\end{equation}
In particular, if $f$ is $\xi_2$-even (odd), then $\Omega$ is $\xi_2$-odd (even).

\end{enumerate}}
\end{Proposition}

Below we consider the inverse of $(t\p_t +\ka-\cL)$, which will be frequently used in the construction of our approximate solutions.
\begin{Lemma}\label{s2:lemma-inversion}
Let $\ka\geq1/2$ and $0< t \leq T\ll \nu^{-1}$. Assume that $F^{\EL}$ and $F^{\NS}$ satisfy
\begin{equation*}
    (t\p_t)^N F^{\EL}, \quad (t\p_t)^N F^{\NS} = \cO_{\cZ}(\log^{M_N}(\nu t)) \quad \text{for any $N\in \N$ and some $M_N\in \N$},
\end{equation*}
and 
\[
 \cP_0F^{\EL}=0,\quad \bar{M}[F^{\EL}](t)=\bar{M}[F^{\NS}](t)=0, \quad  \bar{m}[ F^{\NS}] (t) =0 \quad \text{for $t\leq T$}.
\]
Then there exist 
$\Omega^{\EL},\Omega^{\NS}$ satisfying 
\begin{equation}\label{eq 2.24}
     (t\p_t)^N \Omega^{\EL}, \quad (t\p_t)^N \Omega^{\NS} = \cO_{\cZ}(\log^{M_N}(\nu t)) \quad \text{for any $N\in \N$ and some ${M_N}\in \N$},
\end{equation}
\begin{equation}\label{eq 2.25}
\bar{M}[\Omega^{\EL}](t)=\bar{M}[\Omega^{\NS}](t)=0, \quad \bar{m}[ \Omega^{\EL}] (t)=\bar{m}[ \Omega^{\NS}] (t)=0 \quad \text{for $t\leq T$},
\end{equation}
and 
\begin{equation}\label{eq 2.26}
 \left(t\pt+\ka-\cL-\nu^{-1}\Lam\right)
          (\Omega^{\EL}+\nu\Omega^{\NS})
       +\nu^{-1}F^{\EL}+F^{\NS}
       =\cO^L_{\cZ}(\nu).
\end{equation}
Moreover, the solutions have finite angular Fourier expansions whenever $F^{\EL}$ and $F^{\NS}$ do. The non-radial
part of $\Omega^{\EL}$ has the parity opposite to that of $F^{\EL}$ in
$\xi_2$.
\end{Lemma}

\begin{proof}
We first construct the inverse of $t\pt+\ka-\cL$. Writing $F=Gf$, we define the right inverse by the Duhamel formula and the Ornstein-Uhlenbeck
representation:
\begin{equation}\label{eq 2.27}
\begin{split}
 \mathscr{R}F(t,\xi)
   &\eqdefa \int_0^\oo e^{-\ka s}e^{s\cL}
                         F(te^{-s},\xi)\dd s\\
  & = G(\xi) \int_0^\oo \int_{\R^2} e^{-\ka s}  G(\zeta) f(te^{-s}, e^{-s/2}\xi + \sqrt{1-e^{-s}} \zeta) d\zeta ds.                         
\end{split}
\end{equation}
Differentiating under the integral  gives 
$(t\partial_t+\kappa-\mathcal L)\mathscr R F=F$.

Now for $F^{\EL}$ and $F^{\NS}$ satisfying our assumptions, we set
\begin{equation}\label{eq 2.28}
\begin{aligned}
 \Omega^{\EL}=\Lam^{-1}F^{\EL}
   -\mathscr{R}\cP_0F^{\NS},\quad 
 \Omega^{\NS}=\Lam^{-1}
             \bigl[(t\pt+\ka-\cL)\Omega^{\EL}+F^{\NS}\bigr].
\end{aligned}
\end{equation}
Equations~\eqref{eq 2.27} and~\eqref{eq 2.28} yield \eqref{eq 2.24}. The identities in \eqref{eq
2.25} follow from Proposition \ref{Prop a1}, $\bar{m}[e^{s\cL}h]=\bar{m}[h]$, $\bar{M}[e^{s\cL}h]=e^{-\f{s}{2}}\bar{M}[h] $. Substitution into the left-hand side of \eqref{eq 2.26} leaves exactly
\[
 \nu(t\pt+\ka-\cL)\Omega^{\NS}
           =\cO^L_{\cZ}(\nu),
\]
which proves \eqref{eq 2.26}.  Finally, $\Lambda^{-1}$ reverses parity in $\xi_2$ on the non-radial subspace. This proves the parity assertion and completes the proof of Lemma~\ref{s2:lemma-inversion}.
\end{proof}

\renewcommand{\theequation}{\thesection.\arabic{equation}}
\setcounter{equation}{0}

\section{Construction of approximate solutions}\label{s3:section}
In this section, for each $M\geq 2$, following the general strategy of \cite{GS2024, DG2026,ZZ26}, we construct approximate solutions $(\Omega_a, \dbr_a, \dbz_a)$  of the form
\begin{equation}\label{eq 3.1}
 \Omega_a^{\M}=\chi(\mu \xi/r_0)\sum_{k=0}^M\mu^k\Omega_k,
 \quad
 \dbr_a^{\M}=\sum_{k=0}^{M-1}\mu^k\dbr_k,
 \quad
 \dbz_a^{\M}=\sum_{k=0}^{M-1}\mu^k\dbz_k, \quad \ve^{\M}_a=\frac{\mu}{\br^{\M}_a},
\end{equation}
and then study the functional relationship underlying Arnold's method.
Here $\chi$ is a radial and smooth cut-off function, which equals to one on $|\xi|\leq1/2$, and zero on
$|\xi|\geq 3/4$,
\begin{equation}\label{eq 3.2}
  \Omega_k=\Omega_k^{\EL}(t,\beta_\mu)+\nu\Omega_k^{\NS}(t,\beta_\mu),\quad
 \dbr_k=\nu\dbr_k^{\NS}(t,\beta_\mu),\quad
 \dbz_k=\dbz_k^{\EL}(t,\beta_\mu)+\nu\dbz_k^{\NS}(t,\beta_\mu),
\end{equation}
with $\Omega^{\EL}_k,\Omega^{\NS}_k=\cO^L_\cZ(1)$, $\dbr^{\NS}_k,\dbz^{\EL}_k,\dbz^{\NS}_k=\cO^L(1)$.
By integrating in time, we recover 
\begin{equation}\label{eq 3.3}
    \br_a^{\M} =r_0 + \sum_{k=0}^{M-1} \mu^{k+2}a_{k+2}(t,\beta_\mu) , \quad \bz_a^{\M}=z_0 + \sum_{k=0}^{M-1} \mu^k b_{k}(t,\beta_\mu),  
\end{equation}
with 
\begin{equation*}
   a_{k+2}=\frac{1}{t\mu^k}\int_0^t\mu(\tau)^k\dot{\bar
r}^{\NS}_k(\tau,\beta_{\mu(\tau)})\,d\tau=\cO^L(1),\qquad b_k=\mu^{-k}\int_0^t\mu(\tau)
^k\dot{\bar z}_k(\tau,\beta_{\mu(\tau)})\,d\tau=\cO^L(t).
\end{equation*}
In what follows, we always drop the superscript $(M)$ on $\ve$.
Due to the dependence on $\beta_\mu$ in \eqref{eq 3.2} and \eqref{eq 3.3},  we modify the coefficient-extraction procedure used in \cite{DG2026,ZZ26}: for any given expansion \eqref{eq 3.1},  we define the following formal series
\begin{equation}\label{eq 3.4}
\begin{split}
    &\ve_{\lam}=\frac{\lam}{\br^{\M}_{a,\lam}},\quad \Omega_{a,\lam}^{\M}=\chi(\lam\xi/r_0)\sum_{k=0}^M\lam^k\Omega_k(t,\beta_\mu),
 \\
 &
 \dbr_{a,\lam}^{\M}=\sum_{k=0}^{M-1}\lam^k\dbr_k(t,\beta_\mu),
 \quad
 \dbz_{a,\lam}^{\M}=\sum_{k=0}^{M-1}\lam^k\dbz_k(t,\beta_\mu),\\
 &\br_{a,\lam}^{\M} =r_0 + \sum_{k=0}^{M-1}\lam^{k+2}a_{k+2}(t,\beta_\mu), \quad \bz_{a,\lam}^{\M}=z_0 + \sum_{k=0}^{M-1} \lam^k b_{k}(t,\beta_\mu).
\end{split}
\end{equation}
To proceed, we define $\Phi^{\M}_a\eqdefa \BS_\ve[\Omega_a^{\M}]$ and the corresponding $\lam$-series 
\begin{equation}\label{eq 3.5}
  \Phi^{\M}_{a,\lam} \eqdefa  \BS^{\M}[\Omega^{\M}_{a,\lam}] \with  \BS^{\M} \eqdefa \sum_{k=0}^\oo \ve_\lam^k \Big(\Big(\log \frac{\br_{a,\lam}^{\M}}{\mu} + L\Big)P_k +Q_k\Big).
\end{equation}
And for a formal series $f_\lam$, we write
\[
 \E_k^\lam[f_\lam] \eqdefa \frac{1}{k!} \frac{d^k f_\lam}{d\lam^k}\Big|_{\lam=0}
       \andf 
 \Pi_M[f_\lam] \eqdefa \sum_{k=0}^M \mu^k \E_{k}^{\lam}[f_\lam].
\]

\begin{Remark}
   In the formal expansions above, the recursively constructed coefficients are held fixed when extracting powers of \(\lambda\). In particular, \(t>0\) and \(\beta_\mu\) are treated as fixed parameters during this operation. The formal operator \(\mathrm{BS}^{\M}\) in \eqref{eq 3.5} is obtained from the kernel expansion of the physical operator \(\mathrm{BS}_{\varepsilon_\lambda}\) by replacing
\[
\log\frac{\bar r_{a,\lambda}^{(M)}}{\lambda}
\quad\text{with}\quad
\log\frac{\bar r_{a,\lambda}^{(M)}}{\mu}, \quad \text{while retaining the geometric factors} \, \varepsilon_\lambda^k
=\left(\frac{\lambda}{\bar r_{a,\lambda}^{(M)}}\right)^k.
\]
Thus, the logarithmic singularity in \(\lambda\) is excluded from the formal coefficient extraction. Although $\BS^{\M}[\Omega]$ is a formal expansion, $\E^\lam_k[\BS^{\M}[\Omega]]$, $\Pi_M[\BS^{\M}[\Omega]]$ are always well-defined and we have $\BS_\ve[\Omega]=\Pi_M[\BS^{\M}[\Omega]] + \cO_{\cS_\star}(\beta_\mu \mu^{M+1})$  by Lemma \ref{lem 2.1}, given $\Omega=\cO_\cZ(1)$ and the parameters $(\ve_a,\br_a^{\M},\bz_a^{\M})$ constructed in Proposition \ref{prop 3.1}. 
\end{Remark}

\subsection{Construction of approximate solutions}
We define
\begin{equation}\label{eq 3.6}
\begin{aligned}
 &\RR[\Omega,\dbr,\dbz]
 \eqdefa (t\pt-\cL)\Omega
       -\ve\p_1\left(\varpi^{-1}\Omega\right)
       -\nu^{-1}\ve\br\dbr\p_1\Omega+\nu^{-1}\left\{\BS_\ve[\Omega]-
           \br\dbz\varpi^2/2,
                         \varpi^{-1}\Omega\right\}\\
 &= (t\pt-\cL)\Omega
       -\mu  \p_1\left(\frac{\Omega}{\br+\mu\xi_1}\right)
       -\nu^{-1}\mu \dbr\p_1\Omega+\nu^{-1}\left\{\BS_{\mu/\br}[\Omega]-
           \frac{\dbz}{2\br}(\br+\mu\xi_1)^2,
                \frac{\br}{\br+\mu\xi_1}\Omega\right\}
\end{aligned}
\end{equation}

\begin{Proposition}\label{prop 3.1}
Fix $s\in(0,1)$ and $M\ge2$. There exist $\nu_1(M,s,r_0)>0$, and approximate solutions in the form \eqref{eq 3.1}-\eqref{eq 3.2}, such that for $0<\nu\le\nu_1$ and $0<t\le\nu^{-s}$:
\begin{enumerate}[label=(a\arabic*),ref=a\arabic*]
\item\label{approx:a1}
 $\RR[\Omega_a^{\M},\dbr_a^{\M},\dbz_a^{\M}]
   =\cO^L_{\cZ}(\nu^{-1}\mu^{M+1}+\nu\mu)$. 
\item\label{approx:a2} $\bar{m}[\Omega_0]=1$ and $\bar{m}[\Omega^{\EL}_k]=\bar{m}[\Omega^{\NS}_k]=0$ for $k\geq1$.
\item\label{approx:a3} $\bar{M}[\Omega^{\EL}_k]=\bar{M}[\Omega^{\NS}_k]=0$ for $k\geq0$.
\item\label{approx:a4} $\Omega_k^{\EL},\Omega_k^{\NS}\in\cZ$, and $\Omega_k^{\EL}$
is even in $\xi_2$ for $k\geq0$.
\item\label{approx:a5} $\Omega_k^{\EL},\Omega_k^{\NS}$ have finite angular Fourier
expansions.
\item\label{approx:a6} $
 (t\p_t)^j\Omega_k^{\EL},(t\p_t)^j\Omega_k^{\NS}=\cO^L_{\cZ}(1) $ and  $ (t\p_t)^j\dbr_k^{\NS},(t\p_t)^j\dbz_k^{\EL},(t\p_t)^j\dbz_k^{\NS}=\cO^L(1),$ for any $j\in \N$.
In particular, $\br_a^{\M}(t)\in[3r_0/4,5r_0/4]$.
\end{enumerate}
Moreover,
\begin{equation}\label{eq 3.7}
\begin{aligned}
 \Omega_0^{\EL}&=G,\quad \Omega_0^{\NS}=0,\quad
 \dbr_0^{\NS}=-\f1{r_0},\quad
 \dbz_0^{\NS}=0,\\
 \dbz_0^{\EL}
  &=\f1{4\pi r_0}
          \left(\log \frac{r_0}{\sqrt{\nu t}}+\f32\log2+
                                    \f{\ga_E-1}{2}\right).
\end{aligned}
\end{equation}
\end{Proposition}

\begin{proof}
We prove the proposition by induction, starting with
\begin{equation}\label{eq 3.8}
 \Omega_0=G,\quad   \Phi^{\scriptscriptstyle(0)}_a = -\Upsilon+\f{\log \frac{r_0}{\sqrt{\nu t}}+\log8-2}{2\pi} + \cO_{\cS_\star}^L(\mu),
\end{equation}
which satisfies (\ref{approx:a1})--(\ref{approx:a6}).

\noindent {\bf Step 1. The first-order approximation.} Let
\[
 \Omega_a^{ \scriptscriptstyle(1)}=\chi(\mu\xi/r_0)(G+\mu\Omega_1),\qquad
 \Omega_1=\Omega_1^{\EL}+\nu\Omega_1^{\NS},\qquad
 \dbr_a^{\scriptscriptstyle (1)}=\nu\dbr_0^{\NS},\qquad
 \dbz_a^{\scriptscriptstyle (1)}=\dbz_0^{\EL}+\nu\dbz_0^{\NS}.
\]
By applying Lemma~\ref{lem 2.1}, we obtain
\begin{equation}\label{eq 3.9}
\begin{aligned}
 \RR[\Omega_a^{\scriptscriptstyle (1)},\dbr_a^{\scriptscriptstyle(1)},\dbz_a^{\scriptscriptstyle(1)}]
 =&\mu\left(t\pt+\f12-\cL-\nu^{-1}\Lam\right)\Omega_1
       +\RR[\Omega_a^{\scriptscriptstyle (0)},\dbr_a^{\scriptscriptstyle(1)},\dbz_a^{\scriptscriptstyle(1)}] +\cO^L_{\cZ}(\nu^{-1}\mu^2),
\end{aligned}
\end{equation}
where
\begin{equation*}
\begin{split}
    &\RR[\Omega_a^{\scriptscriptstyle (0)},\dbr_a^{\scriptscriptstyle(1)},\dbz_a^{\scriptscriptstyle(1)}] =\nu^{-1}\mu F^{\EL} +\mu F^{\NS}  +\cO^L_{\cZ}(\nu^{-1}\mu^2),\\
    &F^{\EL} \eqdefa 1/r_0
       \Big(\Big\{((\beta_{ \mu/r_0}+L)P_1+Q_1)[G]
                  -r_0\dbz_0\xi_1,G\Big\} +\{\Upsilon,\xi_1G\}\Big) \quad \text{is $\xi_2$-odd},\\
  &F^{\NS}\eqdefa -(1/r_0+\dbr^{\NS}_0)\p_1G.
\end{split}
\end{equation*}
Thus by choosing  $\dbr_0^{\NS}=-1/r_0$,
$\dbz_0^{\NS}=0$ and
\begin{equation*}
\begin{split}
 \br_0 \dbz_0^{\EL}
   &=\int_{\R^2}G \p_1\big(((\beta_{ \mu/r_0}+L)P_1+Q_1)[G]\big) \dd\xi
                    +\int_{\R^2}\xi_1\p_1\Upsilon G\dd\xi\\
   &=  \f1{4\pi}
          \left(\log\frac{r_0}{\sqrt{\nu t}}+\f32\log2+
                                    \f{\ga_E-1}{2}\right),
\end{split}
\end{equation*}
we have 
\begin{equation*}
F^{\NS}=0, \quad \cP_0 F^{\EL}=0, \quad \bar{M}[F^{\EL}]=0.
\end{equation*}
Then by choosing $\ka=\f12$ and using Lemma \ref{s2:lemma-inversion},  there exist $\Omega_1^{\EL}$ and  $\Omega_1^{\NS}$, so that (\ref{approx:a1})--(\ref{approx:a6}) are satisfied.

\noindent {\bf Step 2. The $(M+1)$-th approximation.}
Assume that the order-$M$ approximation has been constructed. 
We seek the next approximation using the same expansion and splitting:
\begin{equation*}
\begin{aligned}
 \Omega_a^{\MP}
   & = \Omega_a^{\M}+\chi(\mu\xi/r_0)\mu^{M+1}\Omega_{M+1},\quad \dbr_a^{\MP}
   =\dbr_a^{\M}+\mu^M\dbr_M,\quad 
 \dbz_a^{\MP}
   =\dbz_a^{\M}+\mu^M\dbz_M.
\end{aligned}
\end{equation*}
Assuming that the next-order coefficients
satisfy (a6), we obtain
\begin{equation}\label{eq 3.10}
    \begin{split}
        \RR[\Omega_a^{\MP},\dbr_a^{\MP},
\dbz_a^{\MP}]&=\mu^{M+1}\Big(t\pt+\frac{M+1}{2}-\cL-\nu^{-1}\Lam\Big)\Omega_{M+1} \\
&\quad +  \RR[\Omega_a^{\M},\dbr_a^{\MP},
\dbz_a^{\MP}] +\cO^L_{\cZ}(\nu^{-1}\mu^{M+2}+\nu\mu),
    \end{split}
\end{equation}
where we have the natural decomposition in terms of $\mu$ and $\nu$ using corresponding $\lam$-series \eqref{eq 3.4},
\begin{equation*}
    \RR[\Omega_a^{\M},\dbr_a^{\MP},
\dbz_a^{\MP}] = \nu^{-1}\mu^{M+1}F^{\EL}_{M+1}+\mu^{M+1}F^{\NS}_{M+1} + \cO^L_\cZ(\nu^{-1}\mu^{M+2}+\nu\mu).
\end{equation*}
More precisely,
 \begin{equation}\label{eq 3.11}
\begin{split}
    F^{\EL}_{M+1}   &\eqdefa \E^\lam_{M+1} \left\{\BS^{\M}[\Omega^{\M,\EL}_{a,\lam}]-
           \frac{\dbz^{\MP,\EL}_{a,\lam}}{2\br^{\MP}_{a,\lam}}(\br^{\MP}_{a,\lam}+\lam\xi_1)^2,
                \frac{\br^{\MP}_{a,\lam}}{\br^{\MP}_{a,\lam}+\lam\xi_1}\Omega^{\M,\EL}_{a,\lam}\right\} \\
        & \text{is $\xi_2$-odd by assumption (\ref{approx:a4}), and thus $\cP_0F^{\EL}_{M+1}=0$.}
\end{split}
\end{equation}
\begin{equation}\label{eq 3.12a}
\begin{aligned}
F^{\NS}_{M+1}\eqdefa \E^\lam_{M+1}\Bigg[
&-\lam\p_1\left(
  \frac{\Omega^{\M,\EL}_{a,\lam}}
       {\br^{\MP}_{a,\lam}+\lam\xi_1}
  +\dbr^{\MP,\NS}_{a,\lam}\Omega^{\M,\EL}_{a,\lam}
  \right)\\
&+\left\{
  \BS^{\M}[\Omega^{\M,\EL}_{a,\lam}]
  -\frac{\dbz^{\MP,\EL}_{a,\lam}}{2\br^{\MP}_{a,\lam}}
   (\br^{\MP}_{a,\lam}+\lam\xi_1)^2,\,
  \frac{\br^{\MP}_{a,\lam}\Omega^{\M,\NS}_{a,\lam}}
       {\br^{\MP}_{a,\lam}+\lam\xi_1}
  \right\}\\
&+\left\{
  \BS^{\M}[\Omega^{\M,\NS}_{a,\lam}]
  -\frac{\dbz^{\MP,\NS}_{a,\lam}}{2\br^{\MP}_{a,\lam}}
   (\br^{\MP}_{a,\lam}+\lam\xi_1)^2,\,
  \frac{\br^{\MP}_{a,\lam}\Omega^{\M,\EL}_{a,\lam}}
       {\br^{\MP}_{a,\lam}+\lam\xi_1}
  \right\}
\Bigg].
\end{aligned}
\end{equation}
Moreover, $F^{\EL}_{M+1},F^{\NS}_{M+1}=\cO^L_{\cZ}(1)$ have finitely many
angular modes.

\noindent {\bf Step 2.1. Choice of $(\dbr_M,\dbz_M)$.}

We set 
\begin{equation}\label{eq 3.13}
 \left\{   \begin{aligned}
      \dbr^{\NS}_M &\eqdefa - \E_{M+1}^\lam \int_{\HH_{\ve_\lam}} \frac{\br^{\M}_{a,\lam}}{\br^{\M}_{a,\lam}+\lam\xi_1}\p_2\BS^{\M}[\Omega^{\M,\EL}_{a,\lam}] \,\Omega^{\M,\NS}_{a,\lam}\dd\xi \\
&\quad - \E_{M+1}^\lam \int_{\HH_{\ve_\lam}} \frac{\br^{\M}_{a,\lam}}{\br^{\M}_{a,\lam}+\lam\xi_1}\p_2\BS^{\M}[\Omega^{\M,\NS}_{a,\lam}] \,\Omega^{\M,\EL}_{a,\lam}\dd\xi- \E_{M}^\lam \int_{\HH_{\ve_\lam}} \frac{\Omega^{\M,\EL}_{a,\lam}}{\br^{\M}_{a,\lam}+\lam\xi_1} \dd\xi,\\  
 \dbz^{\EL}_M &\eqdefa \E_{M+1}^\lam \int_{\HH_{\ve_\lam}} \frac{\br^{\M}_{a,\lam}}{\br^{\M}_{a,\lam}+\lam\xi_1}\p_1\BS^{\M}[\Omega^{\M,\EL}_{a,\lam}] \,\Omega^{\M,\EL}_{a,\lam}\dd\xi ,\\
 \dbz^{\NS}_M &\eqdefa \E_{M+1}^\lam \int_{\HH_{\ve_\lam}} \frac{\br^{\M}_{a,\lam}}{\br^{\M}_{a,\lam}+\lam\xi_1}\p_1\BS^{\M}[\Omega^{\M,\EL}_{a,\lam}] \,\Omega^{\M,\NS}_{a,\lam}\dd\xi \\
        &\quad + \E_{M+1}^\lam \int_{\HH_{\ve_\lam}} \frac{\br^{\M}_{a,\lam}}{\br^{\M}_{a,\lam}+\lam\xi_1}\p_1\BS^{\M}[\Omega^{\M,\NS}_{a,\lam}] \,\Omega^{\M,\EL}_{a,\lam}\dd\xi.
    \end{aligned}
\right.
\end{equation}
Then together with assumption (a1)-(a6) and \eqref{eq 3.11} and \eqref{eq 3.12a}, the above definition of $(\dbr^{\NS}_M,\dbz^{\EL}_M,\dbz^{\NS}_M)$ determines
\begin{equation}\label{eq 3.15}
   \Bar{m}[F^{\NS}_{M+1}]=0,\quad  \bar{M}[F^{\EL}_{M+1}]=\Bar{M}[F^{\NS}_{M+1}]=0.
\end{equation}
Moreover, in view of \eqref{eq 3.13}, we get $\dbr_{M}^{\NS},\dbz^{\EL}_M,\dbz^{\NS}_M=\cO^L(1),$ and thus for $\nu$ small enough, $\br^{\M}_a(t)$ remains in $[\frac{3}{4}r_0,\frac{5}{4}r_0]$.

\noindent {\bf Step 2.2. Choice of $\Omega_{M+1}$.}
Now thanks to \eqref{eq 3.11}, \eqref{eq 3.15} and Lemma \ref{s2:lemma-inversion}, there exist $\Omega^{\EL}_{M+1},\Omega^{\NS}_{M+1}$ satisfying  (\ref{approx:a2})--(\ref{approx:a6}), so that
\[
 \RR[\Omega_a^{\MP},\dbr_a^{\MP},
                                      \dbz_a^{\MP}]
  =\cO^L_{\cZ}(\nu^{-1}\mu^{M+2}+\nu\mu),
\]
which proves  (\ref{approx:a1}).

This  completes the induction and thus the proof of Proposition \ref{prop 3.1}.
\end{proof}

\subsection{Functional relationship}\label{s3:functional-subsection}
We study the functional relationship between
\begin{equation}\label{eq 3.16}
 \Psi_a^{\M}=\Phi_a^{\M}
       -\f{\br_a^{\M}\dbz_a^{\M}}{2}(\varpi^2-1),
 \qquad
 Q_a^{\M}=\f{\Omega_a^{\M}}{1+\ve_a\xi_1}.
\end{equation}
By using $\lam$-series \eqref{eq 3.4} and Lemma \ref{lem 2.1}, we expand
\begin{equation}\label{eq 3.17}
    \Psi^{\M}_a=\sum_{k=0}^M \mu^k \Psi_k+\cO^L_{\cS_\star}(\mu^{M+1}),\quad  Q^{\M}_a=\sum_{k=0}^M \mu^k Q_k + \cO^{L}_\cZ(\mu^{M+1}),
\end{equation}
where $\Psi_k \eqdefa \E^\lam_k[\Psi^{\M}_{a,\lam}]$ and $Q_k \eqdefa \E^\lam_k[Q^{\M}_{a,\lam}]$, with 
\begin{equation*}
    \begin{split}
  \Psi^{\M}_{a,\lam} \eqdefa       \Phi_{a,\lam}^{\M}
       -\f{\br_{a,\lam}^{\M}\dbz_{a,\lam}^{\M}}{2}\big((1+\ve_\lam\xi_1)^2-1\big),\quad Q^{\M}_{a,\lam} \eqdefa \Omega^{\M}_{a,\lam}/(1+\ve_\lam \xi_1).
    \end{split}
\end{equation*}
We remark that the definitions of $\Psi_k$ and $Q_k$ are independent of the truncation order $M$ whenever $M\ge
k$.

\begin{Proposition}\label{prop 3.2}
Let $M\geq2$. There exist
$F_0(t,\beta_\mu,\cdot),\ldots,F_M(t,\beta_\mu,\cdot)$ satisfying 
\begin{equation}\label{eq 3.18}
    \begin{split}
        (t\p_t)^jF_k(t,\beta_\mu,\cdot)= \cO^L_\cK(1) ,\quad \text{for $j\in\N, k=1,\cdots,M$,}
    \end{split}
\end{equation}
so that
$F^{\M}=\sum_{k=0}^M\mu^kF_k(t,\beta_\mu,\cdot)$ satisfies
\begin{equation}\label{eq 3.19}
 \Pi_M\bigl[\Psi_{a,\lam}^{\M}-F^{\M}_\lam (Q_{a,\lam}^{\M})\bigr]
               =\cO^L_{\cS_\star}(\nu\mu), \with F^{\M}_\lam\eqdefa \sum_{k=0}^M\lam^kF_k(t,\beta_\mu,\cdot).
\end{equation}
\end{Proposition}
\begin{Remark}
    As remarked in \cite{DG2026, ZZ26}, in what follows we only define $F^{\M} (Q_{a}^{\M})$ on $\{|\xi|\lesssim \ve^{-\sig_1}\}$ for $\sig_1\ll1$, wherein $|Q^{\M}_a-G|\ll G$. However,  $\Pi_M\bigl[F^{\M}_\lam (Q_{a,\lam}^{\M})\bigr]$ are well-defined on $\R^2$ since it only involves terms like $\frac{d^k F_k}{ds^k}(Q_0)(Q^{\M}_a-Q_0)^k$.
\end{Remark}
\begin{proof}
Recalling that $Q_0=G$,
\begin{equation}\label{eq 3.20}
\begin{split}
    \Psi_0 &= -\Upsilon+\f{\log \frac{r_0}{\sqrt{\nu t}}+\log8-2}{2\pi},
\end{split}
\end{equation}
we begin by defining
\begin{equation}\label{eq 3.21}
 F_0(t,\beta_\mu,s)=\f1{4\pi}\left(
       2\log r_0 + 2\beta_\mu +4\log2
        -4+\ga_E
        -\Ein\!\left(\log\f1{4\pi s}\right)\right),
                 \qquad 0<s\leq\f1{4\pi}.
\end{equation}
Then $\Psi_0-F_0(Q_0)=0$. To simplify the notation, we set
\begin{equation}\label{eq 3.22}
\begin{split}
   &\qquad \qquad \Psi \eqdefa \sum_{k=0}^M \mu^k \Psi_k, \quad Q \eqdefa \sum_{k=0}^M \mu^k Q_k,  \quad  F=\sum_{k=0}^M\mu^k F_k, \\
   & \Tilde{\Psi} \eqdefa \Psi + \mu^{M+1}\Psi_{M+1}, \quad \Tilde{Q}\eqdefa Q +  \mu^{M+1}Q_{M+1}, \quad \Tilde{F} \eqdefa F+\mu^{M+1}F_{M+1},
\end{split}
\end{equation}
and the subscript $\lambda$ denotes the corresponding formal series obtained by replacing the
explicit powers of $\mu$ with powers of $\lambda$, while keeping all coefficient functions fixed.

We prove by induction. 
Assume that  $\Pi_M[\Psi_\lam-F_\lam(Q_\lam)]=\cO_{\cS_\star}^L(\nu\mu)$. To pass to order $M+1$, we write
\begin{equation*}
    \Pi_{M+1}[\Psi_\lam-F_\lam(Q_\lam)] = \mu^{M+1}\cH_{M+1} +\cO_{\cS_\star}^L(\nu \mu) \quad \text{for some $\cH_{M+1}\in \cS_\star$},
\end{equation*}
which implies that 
\begin{equation}\label{eq 3.23}
    \Pi_{M+1}\{\Psi_\lam,Q_\lam\} = \mu^{M+1} \{\cH_{M+1}, Q_0\} + \cO^L_\cZ(\nu \mu).
\end{equation}
As a result, it follows from  the strategy of constructing $(\Omega_{M+1},\dbr_{M},\dbz_{M})$ that 
\begin{equation*}
\begin{split}
 \cO_\cZ^L(\nu) &= \E^\lam_{M+1}\{\Tilde{\Psi}_\lam,\Tilde{Q}_\lam\} = \{\cH_{M+1}, Q_0\} + \{\Psi_{M+1},Q_0\}+\{\Psi_0,Q_{M+1}\} + \cO_\cZ^L(\nu)\\
 &=\{\cH_{M+1}+\Psi_{M+1}-F_0'(Q_0)Q_{M+1}, Q_0\}  + \cO_\cZ^L(\nu),
\end{split}
\end{equation*}
which together with Proposition \ref{prop 3.1}, implies that  there exists $F_{M+1}$ satisfying \eqref{eq 3.18} so that 
\begin{equation}\label{eq 3.24}
\begin{split}
  &\cP_0\big(\cH_{M+1}+\Psi_{M+1}-F_0'(Q_0)Q_{M+1}\big) =  F_{M+1}(Q_0), \\
  &(1-\cP_0) \big(\cH_{M+1}+\Psi_{M+1}-F_0'(Q_0)Q_{M+1}\big)=\cO_{\cS_\star}^L(\nu).
\end{split}
\end{equation}
Thus by setting $\Tilde{F}=F+\mu^{M+1}F_{M+1}$,  we find 
\begin{equation}\label{eq 3.25}
\begin{split}
    \Pi_{M+1}[\Tilde{\Psi}_\lam-\Tilde{F}_\lam(\Tilde{Q}_\lam)] &= \Pi_{M+1}[\Psi_\lam-F_\lam(Q_\lam)] + \Pi_{M+1}[\Tilde{\Psi}_\lam-\Psi_\lam] - \Pi_{M+1}[\Tilde{F}_\lam(\Tilde{Q}_\lam)-F_\lam(Q_\lam)]\\
    &=\mu^{M+1}({\cH_{M+1}+\Psi_{M+1}-F_0'(Q_0)Q_{M+1}-F_{M+1}(Q_0)})+ \cO_{\cS_\star}^L(\nu\mu)\\
    &=\cO_{\cS_\star}^L(\nu\mu),
\end{split}
\end{equation}
which completes the proof of Proposition \ref{prop 3.2}.
\end{proof}

In what follows, we will use the small parameter $\ve$ instead of $\mu$, the two being comparable.

\begin{Proposition}\label{prop 3.3}
Let $M\geq2$. There exist $\sig_1\ll1$, $N\in\N$, and $C>0$, so that, for sufficiently small
$\ve$ and $|\xi|\leq4\ve^{-\sig_1}$,
\begin{align}
 &|Q_a^{\M}-G(\xi)|
   \leq C\ve\beta_\ve^N\bxi^NG(\xi),\label{eq 3.26}\\
 &|(F^{\M})'(Q_a^{\M})-W_0(\xi)|
   \leq C\ve\beta_\ve^N\bxi^NW_0(\xi),\label{eq 3.27}\\
 &\Theta_a(\ve,\xi)\chi(\ve^{\sig_1}\xi/4)
   =\cO^L_{\cS_\star}
            (\ve^{M+1}+\nu\ve),\label{eq 3.28}
\end{align}
where
\[
 \Theta_a=\Psi_a^{\M}-F^{\M}(Q_a^{\M}),\qquad
 W_0(\xi)\eqdefa F_0'(Q_0)=4|\xi|^{-2}(e^{|\xi|^2/4}-1).
\]
\end{Proposition}

\begin{proof}
\eqref{eq 3.26} and \eqref{eq 3.27} are direct consequences of Propositions \ref{prop 3.1} and \ref{prop 3.2}, by noticing that $W_0(\xi)=F'_0(Q_0)$ and then following the same lines in \cite{GS2024,DG2026,ZZ26}.

And by truncating the expansion \eqref{eq 3.17} at order $M$, and then applying the Taylor formula to the resulting finite-order composition, 
we get from Lemma \ref{lem 2.1}, Propositions \ref{prop 3.1}, \ref{prop 3.2} and $\mu\approx \ve$ that 
\begin{equation*}
\begin{split}
    \chi(\ve^{\sig_1}\xi/4)\Theta_a(\ve,\xi) &= \chi(\ve^{\sig_1}\xi/4) \Pi_{M}\Big[ \Psi_{a,\lam}^{\M} - F^{\M}_\lam(Q^{\M}_{a,\lam})\Big] + \chi(\ve^{\sig_1}\xi/4) \Big(\Psi^{\M}_a - \Pi_M\big[ \Psi^{\M}_{a,\lam}\big]\Big)  \\
    &\quad - \chi(\ve^{\sig_1}\xi/4)\Big(F^{\M}(Q^{\M}_a)-\Pi_M\big[F^{\M}_\lam(Q^{\M}_{a,\lam})\big] \Big) \\
    &= \cO^L_{\cS_\star}
            (\ve^{M+1}+\nu\ve).
\end{split}
\end{equation*}
 We thus complete the proof of Proposition \ref{prop 3.3}.
\end{proof}

\renewcommand{\theequation}{\thesection.\arabic{equation}}
\setcounter{equation}{0}

\section{Estimation of the error}\label{s4:section}

\subsection{Equation for the error}\label{en19:q-section}
We henceforth omit the superscript $(M)$ on the approximate solutions.
We decompose the modulation parameters and define the vorticity and
stream-function errors by
\begin{equation}\label{eq 4.1}
 q\eqdefa\f{\Omega}{1+\ve\xi_1}-Q_a,\qquad
 \br(t)=\br_a(t),\qquad
 \bz_p\eqdefa\bz-\bz_a,\qquad
 \phi\eqdefa \BS_\ve[\varpi q].
\end{equation}
Then we find
\begin{equation}\label{eq 4.2}
\begin{split}
 &(t\pt-t\dbr_a/\br_a
       -\cL-3\ve\varpi^{-1}\p_1)q
 -\nu^{-1}\ve\br_a\dbr_a\p_1q
 +\nu^{-1}\Lam_a[q]
 -\nu^{-1} \ve\br_a\dbz_p\p_2Q_a\\
 &\hspace{25mm}=-\nu^{-1}\varpi^{-1}\{\phi,q\}
 +\nu^{-1}\ve\br_a\dbz_p\p_2q-\RR_a,
\end{split}
\end{equation}
where
\begin{equation}\label{eq 4.3}
 \Lam_a[q]\eqdefa\varpi^{-1}\bigl(\{\Psi_a,q\}+\{\phi,Q_a\}\bigr),
 \qquad
 \RR_a\eqdefa\varpi^{-1}\RR[
       \Omega_a,\dbr_a,\dbz_a].
\end{equation}
Due to the smooth axisymmetric structure of $\omega_\theta/r$ at positive times,  the boundary conditions are
\begin{equation}\label{eq 4.4}
 \p_1q=0,\qquad \phi=\p_1\phi=0
 \qquad\text{on }\p\HH_\ve.
\end{equation}
Here $\bz_p$ is chosen so that $M_2[\Omega]=0$ and $\bz_p(0)=0$, which is equivalent to 
\begin{equation}\label{eq 4.5a}
    \bar z(t)=
\frac{\int_{\HH}z\,\omega_\theta(t,r,z)\,drdz}
{\int_{\HH}\omega_\theta(t,r,z)\,drdz}.
\end{equation}

\subsection{Energy functional and weight function}\label{w19:subsection}
Let
\begin{equation}\label{eq 4.5}
 0<\sig_1\ll 1\ll\sig_2,
 \qquad \ga\eqdefa\f{\sig_1}{\sig_2}\ll1.
\end{equation}
After fixing the approximation order, choose $\sigma_1>0$ sufficiently small that $1-N\sigma_1>\f12$ 
for all polynomial growth exponents used below. Logarithmic factors are then absorbed by a small 
loss in a positive power of $\varepsilon$.
We split $\HH_\ve$ into the three regions
\begin{align}
 \ive&\eqdefa
\{\xi\in\HH_\ve:|\xi|<2\ve^{-\sig_1},\quad
 F'(Q_a(\xi))<\exp(\ve^{-2\sig_1}/4)\},
 \label{eq 4.6}\\
 \iive&\eqdefa\{\xi\in\HH_\ve\setminus \ive:
                  \varrho_\ve\le\ve^{-\sig_2}\},
 \label{eq 4.7}\\
 \iiive&\eqdefa\{\xi\in\HH_\ve:
                  \varrho_\ve>\ve^{-\sig_2}\}.
 \label{eq 4.8}
\end{align}
Here
\[
 \varrho_\ve\eqdefa
 \bigl(\xi_2^2+\chi_0(1+\ve\xi_1)\xi_1^2
       +(1-\chi_0(1+\ve\xi_1))\ve^{-2}\bigr)^{1/2} \approx |\xi|,
\]
where $\chi_0$ is a smooth cut-off function satisfying
$\chi_0(y)=1$ for $y\ge1/2$, and $\chi_0(y)=0$ for $y\le1/4$.
We set
\begin{equation}\label{eq 4.9}
 W_\ve\eqdefa\varpi(\xi)H_\ve,\qquad
 H_\ve\eqdefa
 \begin{cases}
 F'(Q_a),&\xi\in \ive,\\
 \exp(\ve^{-2\sig_1}/4),&\xi\in \iive,\\
 \exp(\varrho_\ve^{2\ga}/4),&\xi\in \iiive.
 \end{cases}
\end{equation}
The weight $H_\ve$ is continuous and locally Lipschitz in space and time, and independent of $\xi_1$ when $0<\varpi(\xi)\le1/4$. We therefore define
\begin{align}
 \|q\|_{X_\ve}^2&\eqdefa
       \int_{\HH_\ve}W_\ve q^2\dd\xi,
 \quad E_\ve[q]\eqdefa\f12\|q\|_{X_\ve}^2
       -\f12\int_{\HH_\ve}\phi q\varpi(\xi)\dd\xi,
 \label{eq 4.10}\\
 \|q\|_{D_\ve}^2&\eqdefa
       \int_{\HH_\ve}W_\ve
        (|\grad q|^2+\rho_\ga^2q^2+q^2)\dd\xi
       +\ve\int_{\p\HH_\ve}H_\ve q^2\dd\xi_2.
 \label{eq 4.11}
\end{align}
Here and below,
$\langle f,g\rangle\eqdefa \int_{\HH_\ve}fg\dd\xi$, and
\begin{equation}\label{eq 4.12}
 \rho_\ga\eqdefa
 \begin{cases}
 |\xi|,&\xi\in \ive,\\
 \ve^{-\sig_1},&\xi\in \iive,\\
 \varrho_\ve^\ga,&\xi\in \iiive.
 \end{cases}
\end{equation}

We summarize some basic properties of the weight $H_\varepsilon$.

\begin{Proposition}[\cite{GS2024}]\label{prop 4.1}
For sufficiently small $\ve$, $\ive$ is diffeomorphic
to an open disk, and there exists $\ka>0$ so that
\begin{equation}\label{eq 4.13}
 \{|\xi|\le\ve^{-\sig_1}\}\subset \ive
 \subset\{|\xi|^2\le\ve^{-2\sig_1}+\ka\beta_\ve\}.
\end{equation}
Moreover, there are universal constants $C_1,C_2,C_3>0$ and $N\in\N$ so that
\begin{align}
 C_1\exp(\varrho_\ve^{2\ga}/4)
   &\le H_\ve\le C_2W_0(\xi), &&\text{for $\xi\in\HH_\ve$},\label{eq 4.14}\\
 |t\p_tH_\ve|+ |H_\ve-W_0|+|\grad H_\ve-\grad W_0|
   &\le C_3\beta_\ve^N\ve^{1-N\sig_1}W_0,  &&\text{for $\xi\in\ive$}.\label{eq 4.15}
\end{align}
\end{Proposition}
\begin{proof}
   In view of the comparability $\varrho_\ve \approx |\xi|$ and the
$t\p_t$-estimates for the approximate solutions and for $F_k$,
$k=1,\ldots,M$, the proof of Proposition \ref{prop 4.1} follows
the argument in \cite{GS2024} with only minor modifications.
We therefore omit the details.
\end{proof}

\begin{Lemma}\label{c19:hardy}
There holds
\begin{equation}\label{eq 4.16}
 \ve\int_{\HH_\ve}
                        |\xi_1H_\ve| q^2\dd\xi
       \lesssim\ve^{-1}\|q\|_{D_\ve}^2.
\end{equation}
\end{Lemma}
\begin{proof}
Since $\ve|\xi_1|\leq 1+ \varpi$,
\[
 \ve\int_{\HH_\ve} |\xi_1H_\ve| q^2\dd\xi
 \le\|q\|_{X_\ve}^2
            +\int_{\HH_\ve}H_\ve q^2\dd\xi,
\]
it suffices to prove
\begin{equation}\label{eq 4.17}
 \ve\int_{\HH_\ve}H_\ve q^2\dd\xi
 \lesssim\int_{\HH_\ve}
                  W_\ve(|\p_1q|^2+q^2)\dd\xi.
\end{equation}
On $\{\varpi\ge1/4\}$ this follows from $H_\ve\le4W_\ve$.
On the remaining strip we use
\begin{equation}\label{eq 4.18}
\int_0^L|g(s)|^2\dd s
 \le C\int_0^L s\bigl(|g'(s)|^2+|g(s)|^2\bigr)\dd s,
 \qquad L\ge2,
\end{equation}
with a constant independent of $L$.
Indeed, integration by parts on $[0,2]$ gives
\[
 \int_0^2g^2\dd s=2g(2)^2-2\int_0^2sgg'\dd s
 \le2g(2)^2+\int_0^2s(g^2+|g'|^2)\dd s.
\]
The elementary trace estimate
$g(2)^2\le C\int_1^2(g^2+|g'|^2)\dd s$ controls the first term.
On $[2,L]$ one simply uses $1\le s/2$. This proves
\eqref{eq 4.18}, first for smooth $g$ and then by
approximation. By setting $L=1/(4\ve)$ and
$g(s)=q(s-\ve^{-1},\xi_2)$, we obtain
\begin{equation}\label{eq 4.19}
 \int_{-\ve^{-1}}^{-3/(4\ve)}
                          |q(\xi_1,\xi_2)|^2\dd\xi_1
 \lesssim\f1\ve
 \int_{-\ve^{-1}}^{-3/(4\ve)}\varpi(\xi)
                  \bigl(|\p_1q|^2+q^2\bigr)(\xi_1,\xi_2)\dd\xi_1.
\end{equation}
Since $H_\varepsilon$ is independent of $\xi_1$ on this strip, multiplying \eqref{eq 4.19} by
$H_\varepsilon$ and integrating in $\xi_2$ yields \eqref{eq 4.17}, and hence \eqref{eq 4.16}.
\end{proof}

\begin{Lemma}\label{lem 4.7}
There holds
\begin{equation}\label{eq 4.20}
\int_{\HH_\ve}|\p_1H_\ve|q^2\dd\xi
 \lesssim\int_{\HH_\ve}W_\ve
                    (|\grad q|^2+\rho_\ga^2q^2+q^2)\dd\xi.
\end{equation}
\end{Lemma}
\begin{proof}
Directly from the definition,
\begin{equation}\label{eq 4.21}
 \p_1H_\ve=0
 \quad  \text{on }
       \{\xi\in\HH_\ve:0<\varpi(\xi)\le1/4\}\cup \iive.
\end{equation}
And Proposition \ref{prop 4.1} implies that, on $\ive$,
\begin{equation}\label{eq 4.22}
 |\p_1H_\ve|
 \lesssim(1+|\xi|)H_\ve
 \lesssim(1+\rho_\ga)W_\ve,
\end{equation}
and on $\iiive\cap\{\varpi>1/4\}$,
\begin{equation}\label{eq 4.23}
 |\p_1H_\ve|
 \lesssim|\p_1\varrho_\ve|
                  \varrho_\ve^{2\ga-1}H_\ve
 \lesssim\ve^{(1-2\ga)\sig_2}W_\ve.
\end{equation}
Finally, integrating the three estimates and
using $1+\rho_\ga\lesssim1+\rho_\ga^2$ proves Lemma \ref{lem 4.7}.
\end{proof}

For convenience, define
\[
 \mathfrak G_\ve[q]\eqdefa W_\ve q-\varpi\phi,
 \qquad \cB_\ve[\cdot]\eqdefa\varpi \BS_\ve[\varpi \cdot].
\]
\begin{Proposition}\label{prop 4.3}
Let $q$ solve \eqref{eq 4.2}--\eqref{eq 4.4}. Then
\begin{equation}\label{eq 4.24}
 (t\p_t-2t\dbr_a/\br_a)E_\ve[q]
   =D_\ve[q]+\cJ_1+\cJ_2+ \cJ_3+\cJ_4,
\end{equation}
where
\begin{align*}
 &D_\ve[q]\eqdefa
  \big\langle(\cL+3\ve\varpi^{-1}\p_1)q,
                     \mathfrak G_\ve[q]\big\rangle
                 +\f12\langle q^2,t\p_tW_\ve\rangle -\f12\langle q,[t\pt,\cB_\ve]q\rangle,\quad \\
 &\cJ_1\eqdefa-\nu^{-1}\langle\Lam_a[q],\mathfrak G_\ve[q]\rangle, \quad 
 \cJ_2\eqdefa-\langle \RR_a,\mathfrak G_\ve[q]\rangle, \\
 &\cJ_3\eqdefa \nu^{-1}\left\langle
   \ve\br_a\dbr_a\p_1q
  +\ve\br_a\dbz_p\p_2Q_a,
                             \mathfrak G_\ve[q]\right\rangle,\\
 &\cJ_4\eqdefa-\nu^{-1}\langle\varpi^{-1}\{\phi,q\},\mathfrak G_\ve[q]\rangle
   +\ve\nu^{-1}\br_a\dbz_p\langle\p_2q,\mathfrak G_\ve[q]\rangle.
\end{align*}
\end{Proposition}
\begin{proof}
Since the kernel of $\cB_\ve$ is symmetric, $W_\varepsilon=\phi=0$ on $\partial\mathbb H_\varepsilon$, differentiating the
quadratic energy gives
\[
 t\p_tE_\ve[q]
 =\langle t\p_tq,\mathfrak G_\ve[q]\rangle
 +\f12\langle q^2,t\p_tW_\ve\rangle
 -\f12\langle q,[t\pt,\cB_\ve]q\rangle.
\]
Then substitution of \eqref{eq 4.2} proves the identity and thus completes the proof of Proposition \ref{prop 4.3}.
\end{proof}

\subsection{Coercivity of $E_\ve$.} We study the coercivity of $E_\ve$ in this subsection.
\begin{Proposition}\label{prop 4.2}
There exist $C_6,C_7>0$ so that, for sufficiently small $\ve$,
\begin{equation}\label{eq 4.25}
 \|q\|_{X_\ve}^2
 \le C_6E_\ve[q]
       +C_7\bigl(\beta_\ve m[q\varpi]^2+|M[q\varpi]|^2\bigr).
\end{equation}
\end{Proposition}
\begin{proof}
Let $\eta=q\varpi$ and define
\begin{equation}\label{eq 4.26}
 E'_\ve[\eta]\eqdefa E_\ve[q]
            =\f12\|\eta\|_{Y_\ve}^2-E^{\mathrm{kin}}[\eta],
\end{equation}
where
\begin{equation}\label{eq 4.27}
 \|\eta\|_{Y_\ve}^2\eqdefa
   \int_{\HH_\ve}\f{H_\ve}{\varpi(\xi)}\eta^2\dd\xi,
 \qquad E^{\mathrm{kin}}[\eta]\eqdefa\f12\langle\phi,\eta\rangle
        =\f12\int_{\HH_\ve}
                      \f{|\grad\phi|^2}{\varpi(\xi)}\dd\xi.
\end{equation}
We use the core--tail decomposition of
\cite[proof of Proposition 4.6]{GS2024}, keeping the singular weight
$H_\ve/\varpi$ in all estimates.
Set
\[
 B_\ve\eqdefa\{\xi\in\HH_\ve:
                           |\xi|\le\ve^{-\sig_1}\},\qquad
 \eta_1\eqdefa\eta\mathds{1}_{B_\ve},\qquad
 \eta_2\eqdefa\eta\mathds{1}_{\HH_\ve\setminus B_\ve}.
\]
Then
\begin{equation}\label{eq 4.28}
 \begin{split}
 E'_\ve[\eta]
 ={}&\f12\|\eta_1\|_{Y_\ve}^2
    +\f12\|\eta_2\|_{Y_\ve}^2
    -E^{\mathrm{kin}}[\eta_1]-E^{\mathrm{kin}}[\eta_2]
    -\int_{\HH_\ve}\phi_1\eta_2\dd\xi,
 \end{split}
\end{equation}
where $\phi_1=\BS_\ve[\eta_1]$. Then in view of \eqref{eq 4.14}, weighted Cauchy-Schwarz inequality and Lemma \ref{lem 2.2}, we get for every $N$,
\begin{align}
 |m[\eta_1]-m[\eta]|^2+|M[\eta_1]-M[\eta]|^2 +|E^{\mathrm{kin}}[\eta_2]|
 +\left|\int_{\HH_\ve}\phi_1\eta_2\dd\xi\right|
     &\le C_N\ve^N\|\eta\|_{Y_\ve}^2.
     \label{eq 4.29}
\end{align}
On $B_\ve$, \eqref{eq 4.15} and
$\varpi^{-1}=1+O(\ve^{1-\sig_1})$ give
\begin{equation}\label{eq 4.30}
 \left|H_\ve/\varpi-W_0\right|
 \le C\bigl(\beta_\ve^N\ve^{1-N\sig_1}
                   +\ve^{1-\sig_1}\bigr)W_0.
\end{equation}
And following a similar idea in \cite{GS2024}, we expand using Lemma \ref{lem 2.1} that 
\begin{equation}\label{eq 4.31}
 \begin{split}
 E^{\mathrm{kin}}[\eta_1]
 ={}&\f{\beta_\ve-2}{4\pi}m[\eta_1]^2
   +\f1{4\pi}\iint_{\R^2\times\R^2}
      \log\f8{|\xi-\zeta|}\eta_1(\xi)\eta_1(\zeta)\dd\xi\dd\zeta\\
 &+O\left(\ve\beta_\ve
                      \int_{\R^2}W_0\eta_1^2\dd\xi\right).
 \end{split}
\end{equation}
Here $\eta_1$ is extended by zero. 
The planar coercivity estimate
\cite[Proposition 4.5]{GS2024} states that, for some $C_1,C_2>0$,
\begin{equation}\label{eq 4.32}
\begin{split}
 \int_{\R^2}W_0\eta_1^2 \dd\xi
 \le{}&C_1\left(\f12\int_{\R^2}W_0\eta_1^2 \dd\xi
 -\f1{4\pi}\iint_{\R^2\times\R^2}
       \log\f8{|\xi-\zeta|}\eta_1(\xi)\eta_1(\zeta)\dd\xi\dd\zeta\right)\\
 &+C_2\bigl(m[\eta_1]^2+|M[\eta_1]|^2\bigr).
\end{split}
\end{equation}
By combining \eqref{eq 4.29}--\eqref{eq 4.32}, we obtain
\begin{align*}
 &\|\eta_1\|_{Y_\ve}^2 \lesssim \int_{\R^2} W_0 \eta_1^2 \dd\xi \\
 &\leq C_1\left(\f12\int_{\R^2}W_0\eta_1^2 \dd\xi
 -\f1{4\pi}\iint_{\R^2\times\R^2}
       \log\f8{|\xi-\zeta|}\eta_1(\xi)\eta_1(\zeta)\dd\xi\dd\zeta\right) +C_2\bigl(m[\eta_1]^2+|M[\eta_1]|^2\bigr)  \\
 &\leq  C_1\left(\f12\|\eta_1\|_{Y_\ve}^2
                              -E^{\mathrm{kin}}[\eta_1]\right)
       +C_2'\bigl(\beta_\ve m[\eta]^2+|M[\eta]|^2\bigr)+ C_3(\ve \beta_\ve + \beta_\ve^N\ve^{1-N\sig_1}+\ve^{1-\sig_1}\big)\|\eta\|^2_{Y_\ve}.
\end{align*}
Finally by inserting the above inequality and \eqref{eq 4.29} into \eqref{eq 4.28}, we achieve \eqref{eq 4.25} and thus complete the proof of Proposition \ref{prop 4.2}.  
\end{proof}

\subsection{Coercivity of $D_\ve[q]$}
We study the coercivity of $D_\ve[q]$ in this subsection. Compared with the dissipation functional in \cite{GS2024}, our functional contains the additional
boundary dissipation term $\mathcal D_{\varepsilon,\mathrm{bry}}$.

\begin{Proposition}[Expansion of the diffusive term]\label{en19:diffusion-prop}
There holds
\begin{equation}\label{eq 4.33}
 D_\ve[q]
    =\cD_{\ve,m}[q]+\cD_{\ve,\mathrm{err}}[q]
                         +\cD_{\ve,\mathrm{bry}}[q],
\end{equation}
where
\begin{equation}\label{eq 4.34}
\begin{split}
 \cD_{\ve,m}[q]\eqdefa{}&
       -\int_{\HH_\ve}W_\ve|\grad_\xi q|^2\dd\xi
       -\int_{\HH_\ve}
          (\grad_\xi\log H_\ve\cdot\grad_\xi q)qW_\ve\dd\xi-\int_{\HH_\ve}V_\ve q^2W_\ve\dd\xi,\\
 \cD_{\ve,\mathrm{err}}[q]\eqdefa{}&
       -\ve\int_{\HH_\ve}(\p_1H_\ve)q^2\dd\xi
       +\f{2\ve t\dbr_a}{\br_a}\zeta_\ve[q]
       -\f{t\dbr_a}{2\br_a}\ve
                         \int_{\HH_\ve}\xi_1H_\ve q^2\dd\xi,\\
 \cD_{\ve,\mathrm{bry}}[q]\eqdefa{}&
       -\ve\int_{\p\HH_\ve}H_\ve q^2\dd\xi_2,
\end{split}
\end{equation}
\begin{equation}\label{eq 4.35}
 \zeta_\ve[q]\eqdefa
 \f14\int_{\HH_\ve}
       \f{\xi_1|\grad_\xi\phi|^2}{\varpi^2}\dd\xi
 +\f12\int_{\HH_\ve}\xi_1\phi q\dd\xi,
\end{equation}
\begin{equation}\label{eq 4.36}
 V_\ve\eqdefa
 \f14\xi\cdot\grad_\xi\log H_\ve
 -\f12t\pt\log H_\ve-\f12
 -\f{\varpi^3}{W_\ve}.
\end{equation}
\end{Proposition}
\begin{proof}
    See Appendix \ref{appA}.
\end{proof}

\begin{Proposition}\label{c19:coercivity}
There exist $\wh\ka,C>0$ so that
\begin{equation}\label{eq 4.37}
\cD_{\ve,m}[q]+\cD_{\ve,\mathrm{bry}}[q]
 \le-\wh\ka\|q\|_{D_\ve}^2
       +C\bigl(m[q\varpi]^2+|M[q\varpi]|^2\bigr).
\end{equation}
\end{Proposition}
\begin{proof}
We follow the localization argument in
\cite{GS2024} and give a proof sketch.

Choose a smooth partition of unity
$1=\chi_3^2+\chi_4^2$, where $\chi_3=1$ for
$|\xi|\le\ve^{-\sig_1}/2$ and $\chi_3=0$ for
$|\xi|\ge\ve^{-\sig_1}$, with
$|\grad\chi_3|+|\grad\chi_4|\le C\ve^{\sig_1}$.
Set $q_3=q\chi_3$ and $q_4=q\chi_4$. Then
\begin{equation}\label{eq 4.38}
 \cD_{\ve,m}[q]
 =\cD_{\ve,m}[q_3]+\cD_{\ve,m}[q_4]
   +\int_{\HH_\ve}W_\ve
          (|\grad\chi_3|^2+|\grad\chi_4|^2)q^2\dd\xi.
\end{equation}
Moreover, the cutoff estimates and the weighted Cauchy--Schwarz inequality give 
\begin{equation}\label{eq 4.39}
 \int_{\HH_\ve}W_\ve
          (|\grad\chi_3|^2+|\grad\chi_4|^2)q^2\dd\xi
 \le C\ve^{2\sig_1}\|q\|_{X_\ve}^2, \andf |m[q_4\varpi]| + |M[q_4\varpi]|\lesssim e^{-c\ve^{-2\sig_1}} \|q\|_{X_\ve}.
\end{equation}
For the core,  by using  Proposition \ref{prop 4.1}, the bound
$|\varpi-1|\leq \ve^{1-\sig_1}$ on the support of $q_3$, and 
\cite[Proposition 4.14]{GS2024}, we obtain
\begin{align*}
 &\|q_3\|_{D_\ve}^2 \approx \int_{\R^2}W_0
       \bigl(|\grad q_3|^2+(1+|\xi|^2)q_3^2\bigr)\dd\xi\notag\\
 &\quad\le C_1\int_{\R^2}
 W_0|\grad q_3|^2+(\grad W_0\cdot\grad q_3)q_3
       +\left(\f14\xi\cdot\grad W_0-\f12W_0-1\right)q_3^2\dd\xi
       +C_2\bigl(m[q_3]^2+|M[q_3]|^2\bigr)\\
&\quad \leq -C_1 \cD_{\ve,m}[q_3] +C_2\bigl(m[q_3\varpi]^2+|M[q_3\varpi]|^2\bigr) + C_3\bigl(\beta_\ve^N\ve^{1-N\sig_1} +\ve^{1-\sig_1}\bigr) \|q_3\|_{D_\ve}^2,
\end{align*}
which implies that 
\begin{equation}\label{eq 4.40}
    \begin{split}
&\|q_3\|_{D_\ve}^2  \leq - C_4 \cD_{\ve,m}[q_3] +C_5\bigl(m[q_3\varpi]^2+|M[q_3\varpi]|^2\bigr) .
    \end{split}
\end{equation}

For the exterior part, Young's inequality gives
\begin{equation}\label{eq 4.41}
 -\cD_{\ve,m}[q_4]
 \ge\f14\|\grad q_4\|_{X_\ve}^2
 +\int_{\HH_\ve}
       \left(V_\ve-\f13|\grad\log H_\ve|^2\right)
                    q_4^2W_\ve\dd\xi.
\end{equation}
By direct calculations, we have 
\begin{equation*}
\begin{split}
    &2\varrho_\ve \p_1 \varrho_\ve= 2\xi_1 \chi_0(\varpi) +\ve \xi_1^2 \chi'_0(\varpi) -\ve^{-1}\chi'_0(\varpi), \\
    & 2\varrho_\ve \p_2 \varrho_\ve =2\xi_2, \\
    & 2\varrho_\ve t\p_t \varrho_\ve = t\pt\ve \xi_1 (\xi_1^2-\ve^{-2}) \chi'_0(\varpi) - 2\ve^{-3} t\pt\ve (1-\chi_0(\varpi))=\cO(\ve^{-2}).
\end{split}
\end{equation*}
Thus for $\sig_2\gg1$
\begin{equation*}
\begin{split}
 V_\ve - \f13 |\grad\log H_\ve|^2 &\geq \begin{cases}
      c_1|\xi|^2 -c_2 , &\quad \xi\in \ive,\\
      \frac{\sig_1}{4}\frac{t\pt\ve}{\ve} \ve^{-2\sig_1} -1/2-\varpi^2/H_\ve, &\quad \xi\in \iive,\\
      \frac{\ga}{8}(\xi\cdot\grad\varrho_\ve - 2t\pt \varrho_\ve) \varrho_\ve^{2\ga-1} - \frac{\ga^2}{12}|\grad\varrho_\ve|^2 \varrho_\ve^{4\ga-2} -1/2-\varpi^2/H_\ve ,&\quad \xi\in \iiive,
  \end{cases}\\
  &\geq  \begin{cases}
  c_1|\xi|^2 -c_2 , &\quad \xi\in \ive,\\
  c_3 \ve^{-2\sig_1} -c_4, &\quad \xi\in \iive,\\
  c_5 \varrho_\ve^{2\ga} -c_6, &\quad \xi\in \iiive,
   \end{cases}  \geq c\rho_\ga^2 -C,
\end{split}
\end{equation*}
which together with \eqref{eq 4.41} and the fact that $\supp q_4 \subset \{|\xi|\geq \ve^{-\sig_1}/2\}$,  implies that 
\begin{equation}\label{eq 4.42}
    \begin{split}
-\cD_{\ve,m}[q_4] &\geq  \f14 \|\grad q_4\|_{X_\ve}^2 + \int_{\HH_\ve}  (V_\ve - \f13 |\grad\log H_\ve|^2) q_4^2 W_\ve d\xi\\
&\geq  \hat{\ka} \Big( \|\grad q_4\|_{X_\ve}^2 + \|\rho_\ga q_4\|_{X_\ve}^2 + \|q_4\|_{X_\ve}^2\Big).
    \end{split}
\end{equation}
By inserting \eqref{eq 4.39}, \eqref{eq 4.40} and \eqref{eq 4.42} into \eqref{eq 4.38}, we obtain
\[
 \cD_{\ve,m}[q]
 \le-\wh\ka\int_{\HH_\ve}W_\ve
              \bigl(|\grad q|^2+(1+\rho_\ga^2)q^2\bigr)\dd\xi
       +C\bigl(m[q\varpi]^2+|M[q\varpi]|^2\bigr).
\]

Finally, add
$\cD_{\ve,\mathrm{bry}}[q]
=-\ve\int_{\p\HH_\ve}H_\ve q^2\dd\xi_2$
and decrease $\wh\ka$ so that $\hat{\ka}\leq 1$ if necessary. This  completes the proof of Proposition \ref{c19:coercivity}.
\end{proof}


The following coercivity inequality is crucial for our stability analysis:
\begin{Proposition}\label{prop 4.9}
There exist $\ka,C>0$ so that, for sufficiently small
$\ve$,
\begin{equation}\label{eq 4.43}
 D_\ve[q]
 \le-\ka\|q\|_{D_\ve}^2
       +C\bigl(m[\varpi q]^2+|M[\varpi q]|^2\bigr).
\end{equation}
\end{Proposition}
\begin{proof}
In view of Proposition \ref{c19:coercivity}, it suffices to prove
$|\cD_{\ve,\mathrm{err}}[q]|\lesssim \ve \|q\|_{D_\ve}^2.$ Indeed, recall that
\[
 \ve\zeta_\ve[q]
 =\f12\int_{\HH_\ve}
                 \eta\,\xi\cdot\grad\phi\dd\xi
  +\f{\ve}{2}\int_{\HH_\ve}
                 \xi_1\phi q\dd\xi,
\]
which together with $\eta=\varpi q$ and Lemma \ref{lem 2.2}, implies that 
\begin{equation*}
\begin{split}
   \ve|\zeta_\ve[q]| &\lesssim  \|\varpi^{-1}\grad\phi\|_{L^4}
           \|\bxi \varpi^2 q\|_{L^{4/3}} + \ve \left\|\f{\bxi\phi}
                         {\sqrt{\varpi H_\ve}}\right\|_{L^2}
            \|q\|_{X_\ve} \\
        &\lesssim \|\bxi^2 \varpi q\|_{L^2} \|\bxi \varpi^2 q\|_{L^{4/3}} + \ve\beta_\ve\|\bxi^2 \varpi q\|_{L^2}\|q\|_{X_\ve}  \lesssim \|q\|_{X_\ve}^2.
\end{split}
\end{equation*}
By combining this bound with Lemmas~\ref{c19:hardy}--\ref{lem 4.7} and using $|t\dot{\bar
r}_a|\lesssim\varepsilon^2$, we obtain
\begin{equation*}
    \begin{split}
|\cD_{\ve,{\rm err}} [q]| \lesssim \ve\|q\|_{D_\ve}^2.
    \end{split}
\end{equation*}
Thus $\cD_{\ve,\mathrm{err}}$ can be absorbed
by the coercive bound for
$\cD_{\ve,\mathrm m}+\cD_{\ve,\mathrm{bry}}$ in
Proposition~\ref{c19:coercivity}. This completes  the proof of Proposition \ref{prop 4.9}.
\end{proof}

\subsection{Estimation of $\cJ$}In this subsection, we estimate $\cJ_1,\ldots,\cJ_4$.

\begin{Lemma}\label{bd19:J1}
For some fixed integer $N$, there holds
\begin{equation}\label{eq 4.44}
 |\cJ_1|\lesssim \nu^{-1}\beta_\ve^N
 \bigl(\ve^{M+1-N\sig_1}
       +\nu\ve^{1-N\sig_1}
       +\ve^{\sig_2}\bigr)\|q\|_{D_\ve}^2.
\end{equation}
\end{Lemma}
\begin{proof}
Recall $\Psi_a=\Phi_a-(\varpi^2-1) \br_a\dbz_a/2$, $\Theta_a=\Psi_a-F(Q_a)$, and $H_\ve =F'(Q_a)$ in $\ive$.
Decompose
\begin{align*}
 \langle\Lam_a[q],\mathfrak G_\ve[q]\rangle
 &=\langle\{\Psi_a,q\}+\{\phi,Q_a\},
                                      qH_\ve-\phi\rangle\\
   &=-\tfrac12\langle\{\Psi_a,H_\ve\},q^2\rangle
       +\langle H_\ve\{\phi,Q_a\}
                -\{\phi,\Psi_a\},q\rangle=\mathscr{A}_{11}+\mathscr{A}_{12}+\mathscr{A}_{21}+\mathscr{A}_{22},
\end{align*}
where
\begin{align*}
 \mathscr{A}_{11}&=-\tfrac12\int_{\ive}
                         \{\Theta_a,H_\ve\} q^2 \dd\xi,
 &\mathscr{A}_{12}&=-\tfrac12\int_{\iiive}
                         \{\Psi_a,H_\ve\} q^2  \dd\xi,\\
 \mathscr{A}_{21}&=-\int_{\ive}\{\phi,\Theta_a\} q  \dd\xi,
 &\mathscr{A}_{22}&=\int_{{\rm I}_{\ve}^c}
    \bigl(H_\ve\{\phi,Q_a\}
                    -\{\phi,\Psi_a\}\bigr)q  \dd\xi.
\end{align*}
Here  we used $\{F(Q_a),F'(Q_a)\}=0$ on $\ive$ and  $\grad_\xi H_\ve$ vanishes on $\iive$.

The functional-relation estimate in Section~\ref{s3:section} gives, after increasing
$N$ if necessary,
\[
 |\grad\Theta_a(\xi)|
 \lesssim \beta_\ve^N
       (\ve^{M+1}+\nu\ve)\bxi^N
       \quad\hbox{on }\ive.
\]
Since $\ive\subset\{|\xi|\le C\ve^{-\sig_1}\}$,
$|\varpi-1|\lesssim \ve^{1-N\sig_1}$ there, and
$|\grad H_\ve|\le C(1+\rho_\ga)H_\ve$, it follows
that
\[
 |\mathscr{A}_{11}|\le C\beta_\ve^N
    (\ve^{M+1-N\sig_1}
             +\nu\ve^{1-N\sig_1})
                             \|q\|_{D_\ve}^2.
\]
For $\mathscr{A}_{21}$, Lemma \ref{lem 2.2} gives
\[
 |\mathscr{A}_{21}|\le\|\varpi^{-1}\grad\phi\|_{L^4}
 \left\|\sqrt{\f{\varpi}{H_\ve}}
        \grad\Theta_a\,\mathds{1}_{\ive}\right\|_{L^4}
                        \|q\|_{X_\ve}
 \le C\beta_\ve^N(\ve^{M+1}+\nu\ve)
                        \|q\|_{X_\ve}^2.
\]
On $\iiive$, by noticing that
$\grad\log H_\ve=(\ga/2)
\varrho_\ve^{2\ga-1}\grad\varrho_\ve$,
$|\grad\varrho_\ve|\le C$ and
$\varrho_\ve^{-1}\le\ve^{\sig_2}$, we obtain
\[
 |\mathscr{A}_{12}|\lesssim  \ve^{\sig_2}\|\varpi^{-1}\grad\Psi_a\|_{L^\oo}  \|\rho_\ga q\|_{X_\ve}^2 \lesssim \beta_\ve^N\ve^{\sig_2}
                             \| q\|_{D_\ve}^2.
\]
Finally, by using Lemma \ref{lem 2.2}, we obtain
\begin{align*}
 |\mathscr{A}_{22}|\le{}&\|\varpi^{-1}\grad\phi\|_{L^4}
                         \|q\|_{X_\ve}\left(
 \|\sqrt{\varpi H_\ve}\grad Q_a\,
                                      \mathds{1}_{{\rm I}_{\ve}^c}\|_{L^4}
 +\left\|\sqrt{\f{\varpi}{H_\ve}}\grad\Psi_a\,
                                      \mathds{1}_{{\rm I}_{\ve}^c}\right\|_{L^4}
 \right)\\
 &\lesssim e^{-c\ve^{-2\sig_1}}\|q\|_{X_\ve}^2.
\end{align*}
 By combining these four bounds, we achieve \eqref{eq 4.44} and thus complete the proof of Lemma \ref{bd19:J1}.
\end{proof}

\begin{Lemma}\label{lem 4.13}
For some fixed integer $N$, there holds
\begin{equation}\label{eq 4.45}
 |\cJ_2|\lesssim \beta_\ve^N
       (\nu^{-1}\ve^{M+1}+\nu\ve)
                                  \|q\|_{X_\ve}.
\end{equation}
\end{Lemma}
\begin{proof}
The residual estimate and the Gaussian decay of its coefficients give
\[
 \|\RR_a\|_{X_\ve}
 \le C\beta_\ve^{N'}
                   (\nu^{-1}\ve^{M+1}+\nu\ve),
\]
which together with Lemma \ref{lem 2.2}, implies
\begin{align*}
 |\langle \RR_a,\mathfrak G_\ve[q]\rangle|
 &\le\|\RR_a\|_{X_\ve}
 \left(\|q\|_{X_\ve}
       +\left\|\sqrt{\f{\varpi}{H_\ve}}\phi\right\|_{L^2}\right)\lesssim \beta_\ve^{N'+1}
    (\nu^{-1}\ve^{M+1}+\nu\ve)
                                  \|q\|_{X_\ve}.
\end{align*}
We thus complete the proof of Lemma \ref{lem 4.13}.
\end{proof}

\begin{Lemma}\label{lem 4.11}
For some fixed integer $N$, there holds
\begin{equation}\label{eq 4.46}
 |\cJ_3|\lesssim \ve\beta_\ve\|q\|_{D_\ve}^2
 +\nu^{-1}\ve
       |\dbz_p|\beta_\ve^N
       (\ve^{M+1}+\nu\ve)\|q\|_{X_\ve}.
\end{equation}
\end{Lemma}
\begin{proof}
By Lemma \ref{lem 2.2} and $|\br_a\dbr_a|\le C\nu$, we get
\begin{align}\label{eq 4.47}
 \left|\left\langle\f{\ve\br_a\dbr_a}{\nu}\p_1q,\mathfrak G_\ve[q]\right\rangle\right|
&\lesssim \ve\|\grad q\|_{X_\ve}
 \left(\|q\|_{X_\ve}
       +\left\|\sqrt{\f{\varpi}{H_\ve}}\phi\right\|_{L^2}\right) \lesssim \ve\beta_\ve\|q\|_{D_\ve}^2.
\end{align}
By using  $[\BS_\ve,\p_2]=0$ and $\p_2\varpi=0$, we obtain
\begin{align*}
 \langle\p_2Q_a,\mathfrak G_\ve[q]\rangle
 &=\int_{\HH_\ve}\varpi q
       \bigl(H_\ve\p_2Q_a-\p_2\Phi_a\bigr)\dd\xi,\with
 \Phi_a=\BS_\ve[\varpi Q_a].
\end{align*}
The frame term in $\Psi_a$ is independent of
$\xi_2$, and hence on $\ive$
\[
 H_\ve\p_2Q_a-\p_2\Phi_a= H_\ve\p_2Q_a-\p_2\Psi_a
   =-\p_2\Theta_a.
\]
Therefore, by the polynomial bound for $\grad\Theta_a$ and the same
Gaussian and exterior estimates used for $\mathscr{A}_{21}$ and $\mathscr{A}_{22}$,
\begin{equation*}
    \begin{split}
&|\langle\p_2Q_a,\mathfrak G_\ve[q]\rangle|\leq \big|\langle q\varpi, \p_2\Theta_a \mathds{1}_{\ive} \rangle\big|+ \big|\langle q\varpi \mathds{1}_{{\rm I}^c_\ve}, H_\ve\p_2Q_a-\p_2\Phi_a \rangle\big|\\
&\lesssim \|q\|_{X_\ve} \big( \big\| (\varpi/H_\ve)^{1/2}\p_2\Theta_a \mathds{1}_{\ive} \big\|_{L^2} +  \big\|(\varpi/H_\ve)^{1/2}
       (H_\ve\p_2Q_a-\p_2\Phi_a) \mathds{1}_{{\rm I}^c_\ve} \big\|_{L^2} \big)\\
    &\lesssim \beta_\ve^N(\ve^{M+1}+\nu\ve)\|q\|_{X_\ve},
    \end{split}
\end{equation*}
by inserting which and \eqref{eq 4.47} into $\cJ_3$, we complete the proof of Lemma \ref{lem 4.11}. 
\end{proof}

\begin{Lemma}\label{lem 4.12}
For some fixed integer $N$, there holds
\begin{equation}\label{eq 4.48}
 |\cJ_4|\lesssim \nu^{-1}\|q\|_{X_\ve}\|q\|_{D_\ve}^2
 +\nu^{-1}\ve|\dbz_p|
                      \|q\|_{X_\ve}\|q\|_{D_\ve}.
\end{equation}
\end{Lemma}
\begin{proof}
To begin with, we claim
\begin{equation}\label{eq 4.49}
    \|\bxi\varpi^{-1}\grad \phi\|_{L^\oo}\lesssim \|q\|_{X_\ve}^\f12\|q\|_{D_\ve}^\f12.
\end{equation}
Indeed,  we write $\grad \eta = \varpi \grad q + \ve q (1,0)^{\rm T}$, to get from \eqref{eq 4.17}
\begin{equation*}
\begin{split}
    \|\bxi^2 \grad \eta\|_{L^2} &\lesssim \| \bxi^2 \varpi \grad q\|_{L^2} + \ve \|\bxi^2 q\|_{L^2}\\
    &\lesssim \| \bxi^2 \varpi \grad q\|_{L^2} + \ve \Big(\int_{\HH_\ve} H_\ve q^2 \dd\xi\Big)^{\f12} \lesssim \|q\|_{D_\ve},
\end{split}
\end{equation*}
which together with  $\|\bxi^2 \eta\|_{L^2}\lesssim \|q\|_{X_\ve}$ and Lemma \ref{lem 2.2}, implies \eqref{eq 4.49}.

Then by using integration by parts, we get
\begin{equation}\label{eq 4.50}
 -\nu^{-1}\langle\varpi^{-1}\{\phi,q\},
                                             \mathfrak G_\ve[q]\rangle
 =\f1{2\nu}\int_{\HH_\ve}
                                  \{\phi,H_\ve\} q^2 \dd\xi,
\end{equation}
which together with $|\grad H_\ve|\le C(1+\rho_\ga)H_\ve$
and \eqref{eq 4.49}, implies
\begin{align*}
 \nu^{-1}|\langle\varpi^{-1}\{\phi,q\},\mathfrak G_\ve[q]\rangle|&\lesssim      \nu^{-1}\|\varpi^{-1}\grad\phi\|_\oo
       \|q\|_{X_\ve}\|q\|_{D_\ve}\\
 &\lesssim \nu^{-1}\|q\|_{X_\ve}^{3/2}
                         \|q\|_{D_\ve}^{3/2}
 \lesssim \nu^{-1}\|q\|_{X_\ve}
                         \|q\|_{D_\ve}^2.
\end{align*}
Since
\begin{equation*}
 \int\varpi\phi\p_2q \dd\xi
 =\int_{\HH_\ve} \BS_\ve[\varpi q]\p_2(\varpi q) \dd\xi=0,
\end{equation*}
one has
\begin{equation}\label{eq 4.51}
 |\langle\p_2q,\mathfrak G_\ve[q]\rangle| = |\langle\p_2q,q W_\ve\rangle|
 \lesssim \|q\|_{X_\ve}\|q\|_{D_\ve}.
\end{equation}
By inserting \eqref{eq 4.50}--\eqref{eq 4.51} into $\cJ_4$, we obtain \eqref{eq 4.48} and thus complete the proof of Lemma \ref{lem 4.12}.
\end{proof}

\renewcommand{\theequation}{\thesection.\arabic{equation}}
\setcounter{equation}{0}

\section{Estimation of the mass and first momentum}\label{s5:section} In view of the coercivity of $E_\ve$, $D_\ve$, we estimate in this section $m[q\varpi]$ and $M[q\varpi]$.

\begin{Lemma}\label{lem 4.14}
Let $T\ll \nu^{-1}$. 
Assume that $\sup_{0<\tau\le T}\|q(\tau)\|_{X_{\ve(\tau)}}\le1$.
Then for $t\leq T$,
\begin{equation}\label{eq 5.1}
 m[q\varpi]=\cO(\nu^{-1}\ve^\oo).
\end{equation}
\end{Lemma}
\begin{proof}
Choose $0\le\chi_2\le1$ with $\chi_2(r)=0$ for $r\le r_0/4$,
$\chi_2(r)=1$ for $r\ge r_0/2$, and
$r_0|\chi_2'|+r_0^2|\chi_2''|\le C$. Define the localized mass by
\begin{equation}\label{eq 5.2}
 J(t)=\int_{\HH}\chi_2(r)\om_\tht(t,r,z)\dd r\dd z
     =\int_{\HH_\ve}\chi_2(\br_a\varpi)\Omega\dd\xi.
\end{equation}
Positivity and the mass identity in Lemma~\ref{lem 1.4} give
$0\le J(t)\le m[\Omega(t)]\le1$ and $J(0)=1$.
By testing \eqref{eq 2.10} with $\chi_2(\br_a\varpi)$ and then using integration by parts, we obtain
\begin{align}\label{eq 5.3}
 tJ'(t)
 &=\int_{\cC_\ve}\left(
       \ve^2\br_a^2\chi_2''(\br_a\varpi)
       -\ve^2\br_a\varpi^{-1}\chi_2'(\br_a\varpi)
       -\f{\ve\br_a}{\nu}\,
          \varpi^{-1}\p_2\Phi\,\chi_2'(\br_a\varpi)
                    \right)\Omega\dd\xi,
\end{align}
where
\[
 \cC_\ve
 =\left\{\xi\in\HH_\ve:
 \left(\f{r_0}{4\br_a}-1\right)\ve^{-1}
 \le\xi_1\le
 \left(\f{r_0}{2\br_a}-1\right)\ve^{-1}\right\} \subset\Big\{-\frac{7}{8\ve}\le\xi_1\le -\frac{1}{4\ve}\Big\}.
\]
By combining the Gaussian localization of $\Omega_a$, the lower bound on $H_\ve$, the hypothesis and    Lemma~\ref{lem 2.2}, we infer that
\begin{equation*}
    \begin{split}
|tJ'(t)|\lesssim \ve^2 \|\Omega \mathds{1}_{\cC_\ve}\|_{L^1} +\nu^{-1}\ve \|\mathds{1}_{\cC_\ve}\p_2 \Phi\|_{L^4}\|\Omega \mathds{1}_{\cC_\ve}\|_{L^{4/3}}\lesssim \nu^{-1} e^{-c\ve^{-2\sig_1}}, 
    \end{split}
\end{equation*}
which implies that 
\begin{equation*}
    tJ'(t) = \cO(\nu^{-1}\ve^{\oo}) \implies 
 J(t)=1+\cO(\nu^{-1}\ve^\oo).
 \end{equation*}
Since $J(t)\le m[\Omega(t)]\le1$ and
$m[\Omega_a]=1+\cO(\ve^\oo)$, subtraction proves
\eqref{eq 5.1}. We thus complete the proof of Lemma \ref{lem 4.14}.
\end{proof}

\begin{Lemma}\label{m19:moments}
Under the assumptions of Lemma~\ref{lem 4.14}, one has for some fixed integer $N\in \N$,
\begin{equation}\label{eq 5.4}
 |M_1[q\varpi]|
 \lesssim \ve\|q\|_{X_\ve}
     +\beta_\ve^N
        (\nu^{-1}\ve^{M+1}+\nu\ve).
\end{equation}
Moreover
\begin{equation}\label{eq 5.5}
 M_2[q\varpi]=-M_2[\Omega_a]=\cO(\ve^\oo),
\end{equation}
and
\begin{equation}\label{eq 5.6}
 |\dbz_p|
 \lesssim      \beta_\ve (\|q\|_{X_\ve}+\|q\|_{X_\ve}^2)
     +\nu e^{-c\ve^{-2\sig_1}} \|q\|_{D_\ve}
     +\beta_\ve^N(\ve^M+\nu^2).
\end{equation}
\end{Lemma}
\begin{proof}
Recall that the axial impulse is conserved:
\[
 \int_{\HH}r^2\om_\tht(t,r,z)\dd r\dd z= r_0^2 \implies \br_a^2(t)
      \int_{\HH_\ve}(\Omega_a+q\varpi)\varpi^2\dd\xi
   = r_0^2.
\]
For the approximate solution, the same computation yields the identity with the residual term
\begin{align*}
 \br_a^2(t)
      \int_{\HH_\ve}\Omega_a\varpi^2\dd\xi
   = r_0^2+
       \int_0^t\f{\br_a^2(\tau)}{\tau}
       \int_{\HH_{\ve(\tau)}}\varpi^3\RR_a\dd\xi\dd\tau,
\end{align*}
which together with $m[\varpi \RR_a]=t\pt m[\Omega_a]=\cO(\ve^{\oo})$ and $\RR_a=\cO_{\cZ}(\nu^{-1}\ve^{M+1}+\nu\ve)$, implies 
\begin{equation}\label{eq 5.7}
\begin{split}
 &\br_a^2(t)\int_{\HH_\ve}q\varpi^3\dd\xi
 =-\int_0^t\f{\br_a^2(\tau)}{\tau}
       \int_{\HH_{\ve(\tau)}}\varpi^3\RR_a\dd\xi\dd\tau \\
& =-\int_0^t\f{\br_a^2(\tau)}{\tau}
       \int_{\HH_{\ve(\tau)}}(\varpi^2-1)\varpi\RR_a\dd\xi\dd\tau -\int_0^t\f{\br_a^2(\tau)}{\tau} m[\varpi \RR_a]\dd\tau\\
    &= \cO^L(\nu^{-1}\ve^{M+2}+\nu \ve^2).
\end{split}
\end{equation}
By using \eqref{eq 5.1} and \eqref{eq 5.7} and expanding
\[
 \int q\varpi^3 \dd\xi
   =m[q\varpi]+2\ve M_1[q\varpi]
                     +\ve^2\int\xi_1^2q\varpi \dd\xi,
\]
we prove \eqref{eq 5.4} under the assumption of Lemma \ref{lem 4.14}. 

To proceed, we define
\begin{equation*}
 \cI[\Omega_1,\Omega_2]
   \eqdefa \int_{\HH_\ve}\varpi^{-1}\p_1\BS_\ve[\Omega_1]\,\Omega_2\dd\xi.
\end{equation*}
By using Lemma \ref{lem 1.4}  and $M_2[\Omega]=0$, we obtain
\begin{equation}\label{eq 5.8}
 (\dbz_a+\dbz_p)m[\Omega_a+q\varpi]
   =\f{1}{\ve\br_a}
                   \cI[\Omega_a+q\varpi,\Omega_a+q\varpi]
     -\f{2\nu}{\br_a}
                   \int_{\R}\xi_2q(-\ve^{-1},\xi_2)\dd\xi_2.
\end{equation}
And a similar computation, together with Proposition~\ref{prop 3.1}, gives
\begin{align}\label{eq 5.9}
 \dbz_a m[\Omega_a]
       -\f{1}{\ve\br_a}\cI[\Omega_a,\Omega_a]
  =\f{\nu}{\ve\br_a}
       \left(M_2[\varpi \RR_a]
           -\left(t\pt+\f12\right)M_2[\Omega_a]\right)=\cO^L(\ve^{M}+\nu^2).
\end{align}
Moreover, by using Lemma \ref{lem 2.2}, we have the following estimate
\begin{equation}\label{eq 5.10}
    |\cI[\Omega_1,\Omega_2]|\lesssim \|\varpi^{-1}\grad \BS_\ve[\Omega_1]\|_{L^4}\|\Omega_2\|_{L^{4/3}} \lesssim \|\bxi^2 \Omega_1\|_{L^2} \|\bxi^2\Omega_2\|_{L^2}.
\end{equation}
And by using Lemmas \ref{lem 2.2}, \ref{lem 2.1}, we get that if ${\rm supp} \Omega_1, {\rm supp} \Omega_2 \subset \{|\xi|\leq 2 \ve^{-\sig_1}\}$, 
\begin{equation}\label{eq 5.11}
\begin{split}
   &\left|\int_{\HH_\ve} (\varpi^{-1} \p_1\BS_\ve[\Omega_1]-\p_1\BS_0[\Omega_1]) \Omega_2 \dd\xi \right| \lesssim \ve\beta_\ve \|\varpi^{-1}\Omega_1\|_{X_\ve}\|\varpi^{-1}\Omega_2\|_{X_\ve}.
\end{split}
\end{equation}
By combining \eqref{eq 5.10}, \eqref{eq 5.11}, the property of $H_\ve$, and 
\begin{equation*}
    \int_{\R^2} \p_1 \BS_0[\Omega_1] \Omega_2 + \p_1 \BS_0[\Omega_2] \Omega_1 \dd \xi =0,
\end{equation*}
we arrive at
\begin{equation*}
    \begin{split}
&|\cI[\Omega_1, \Omega_2]+\cI[\Omega_2, \Omega_1]|\\
&\lesssim  \big|\cI[\chi(\ve^{\sig_1}\xi) \Omega_1, \chi(\ve^{\sig_1}\xi)\Omega_2]+\cI[\chi(\ve^{\sig_1}\xi)\Omega_2, \chi(\ve^{\sig_1}\xi)\Omega_1]\big| + e^{-c\ve^{-2\sig_1}} \|\varpi^{-1}\Omega_1\|_{X_\ve}\|\varpi^{-1}\Omega_2\|_{X_\ve} \\
&\lesssim \ve \beta_\ve \|\varpi^{-1}\Omega_1\|_{X_\ve} \|\varpi^{-1}\Omega_2\|_{X_\ve},
    \end{split}
\end{equation*}
which implies that 
\begin{equation}\label{eq 5.12}
    |\cI[\Omega_a+q\varpi, \Omega_a +q\varpi] - \cI[\Omega_a,\Omega_a]| \lesssim \ve \beta_\ve (\|q\|_{X_\ve}+\|q\|_{X_\ve}^2).
\end{equation}
And from the definition of $\|q\|_{D_\ve}$, we have
\begin{equation}\label{eq 5.13}
 \left|\int_{\R}\xi_2q(-\ve^{-1},\xi_2)\dd\xi_2\right|
 \lesssim e^{-c\ve^{-2\sig_1}}\|q\|_{D_\ve}.
\end{equation}

Finally, by inserting \eqref{eq 5.13}, \eqref{eq 5.12} and \eqref{eq 5.9} into \eqref{eq 5.8}, and then using $m[ q\varpi]=\cO(\nu^{-1}\ve^\oo)$, $|\dbz_a|\lesssim \beta_\ve$ and $m[\Omega_a]=1+\cO(\ve^{\oo})$, we achieve  \eqref{eq 5.6} and thus complete the proof of Lemma \ref{m19:moments}.
\end{proof}

\renewcommand{\theequation}{\thesection.\arabic{equation}}
\setcounter{equation}{0}

\section{Proof of Theorem~\ref{s1:main}}\label{s6:section}

\begin{proof}[Proof of Theorem~\ref{s1:main}]

Fix $s\in(0,1)$, choose 
$\sig_2=M>4/(1-s)$. Once the finite logarithmic exponents
and spatial growth orders in the estimates of Sections~\ref{s3:section}-\ref{s5:section} have been fixed,
choose $\sig_1>0$ sufficiently small that $N\sig_1<1/2$ in all
those estimates.  Throughout the proof, $N$ may vary from line to line, and is  independent of $\nu$ and $t$.
Proposition~\ref{prop 3.1} yields for $0<t\le\nu^{-s}$,
\begin{equation}\label{eq 6.1}
    \br_a \approx r_0,\qquad
 \ve = \sqrt{\nu t}/\br_a \lesssim \nu^{\frac{1-s}{2}},\qquad
 \nu^{-1}\ve^{M+1}\le \nu\ve,
\end{equation}
provided $\nu$ is sufficiently small depending on $s$. Set
\begin{equation}\label{eq 6.2}
 T_m=\sup\left\{t\le\nu^{-s}:
       \|q(\tau)\|_{X_{\ve(\tau)}}
       \le\nu\ve(\tau)
                    \beta_{\ve(\tau)}^{N_1}
       \text{ for every }0<\tau\le t\right\},
\end{equation}
where $N_1\gg1$ will be chosen after $N$ and we finally decrease $\nu_0$.

\noindent {\bf Step 1. Initialization at fixed $\nu$.}
Introduce $\ve_0=\mu/r_0$ and $d(t)=\mu^{-1}(\br_a-r_0,\bz-z_0)$, and define
\begin{equation}\label{eq 6.3}
    f(t,y) \eqdefa \mu^2 \omega_\theta(t,r_0+\mu y_1, z_0 + \mu y_2)= \Omega(t,y-d(t))=(\Omega_a+\eta)(t,y-d(t)), \quad  y_1 >-\ve_0^{-1}.
\end{equation}
Extend $f$ by zero to the rest of $\R^2$ and write $\|h\|_{\cY}^2\eqdefa 4\pi \int_{\R^2} |h|^2 e^{|y|^2/4} \dd y$. Then \cite[Section~4.4, especially the proof of estimate (1.7) in Theorem 1.1]{GS2019} gives for $t\leq t_0\ll1$,
\begin{equation}\label{eq 6.4}
    \|f(t)-G\|_{\cY} \leq C_\nu \ve_0|\log \ve_0|.
\end{equation}
Moreover, \cite[(4.4)]{GS2019} gives
\begin{equation}\label{eq 6.5}
    \|\grad_y f(t)\|_{L^\oo}\leq C_\nu.
\end{equation}
Now in view of $M_2[\Omega]=0$, we have 
\begin{equation*}
    \begin{split}
0=M_2[\Omega] = \mu^2 \int_{\HH_\ve} \xi_2 \omega_\theta(t,\br_a(t)+\mu\xi_1,\bz(t)+\mu\xi_2) \dd\xi = \int_{\R^2} y_2 f(t,y_1, y_2) \dd y + \frac{z_0-\bz(t)}{\mu} \bar{m}[f](t),
    \end{split}
\end{equation*}
which together with $\bar{m}[f]\rightarrow 1$ as $t\rightarrow 0$, \eqref{eq 6.4} and Proposition \ref{prop 3.1}, implies that for $t\leq t_0\ll1$,
\begin{equation}\label{eq 6.6}
    |d(t)|\leq C_\nu \ve \beta_\ve^{N} \quad \text{for some fixed integer $N$.}
\end{equation}

Now we split 
\begin{equation}\label{eq 6.7}
    \|q\|_{X_\ve}^2 = \int_{\{\varpi \geq 1/4\}} \eta^2 \varpi^{-1}H_\ve \dd\xi + \int_{\{\varpi \leq 1/4\}} \eta^2 \varpi^{-1}H_\ve \dd\xi \eqdefa \cQ_1 + \cQ_2.
\end{equation}
For $\cQ_1$, we use \eqref{eq 6.4} and  \eqref{eq 6.6} and $H_\ve \lesssim e^{|\xi+d|^2/4}$, to get
\begin{equation}\label{eq 6.8}
\begin{split}
    \cQ_1 &\lesssim \|f(t,\cdot+d(t))-G\|_{L^2(H_\ve d\xi)}^2+\|\Omega_a-G\|_{L^2(H_\ve d\xi)}^2 \\
    &\lesssim \|f-G\|_{\cY}^2 + |d(t)|^2 + \|\Omega_a-G\|_{L^2(H_\ve d\xi)}^2     \leq C_\nu \ve^2 \beta_\ve^{N}.
\end{split}
\end{equation}
While for $\cQ_2$, we first get from \eqref{eq 6.5} and the properties of our approximate solutions that 
\begin{equation}\label{eq 6.9}
    \|\varpi^{-1}\eta\|_{L^\oo}\leq \ve^{-1}\|\grad_\xi \eta\|_{L^\oo} \leq \ve^{-1}(\|\grad_y f\|_{L^\oo}+ \|\grad_\xi \Omega_a\|_{L^\oo}) \lesssim C_\nu \ve^{-1}.
\end{equation}
We then bound from \eqref{eq 6.4} and \eqref{eq 6.6} that 
\begin{equation}\label{eq 6.10}
    \begin{split}
\cQ_2 &\lesssim \|\varpi^{-1} \eta\|_{L^\oo} \Big(\int_{\{\varpi \leq \f14\}} \eta^2 e^{|\xi|^2/8} \dd\xi\Big)^\f12 \Big( \int_{\{\varpi \leq \f14\}} e^{-|\xi|^2/8} H_\ve^2 \dd\xi \Big)^\f12 \\
&\leq C_\nu \ve^{-1} e^{-c\ve^{-2}}  \Big(\int_{\{-1\leq \ve\xi_1\leq -\f34\}} |f(t,\xi+d(t))-\Omega_a(\xi)|^2  e^{|\xi|^2/8} \dd\xi\Big)^\f12\\
 &\leq C_\nu \ve^{-1} e^{-c\ve^{-2}} \Big( \|\Omega_a-G\|_{\cY} + \|f-G\|_{\cY} + |d(t)| \Big) \leq  C_\nu e^{-c\ve^{-2}} \beta_\ve^{N}.
    \end{split}
\end{equation}
Then by combining \eqref{eq 6.8} and \eqref{eq 6.10}, we have
\begin{equation}\label{eq 6.11}
    \|q\|_{X_\ve}^2 =\cQ_1 +\cQ_2 \leq C_\nu \ve^2\beta_\ve^{N}.
\end{equation}
Thus by choosing $N_1$ sufficiently large, we have $T_m >0$.

\noindent {\bf Step 2. Bootstrap estimates.}
For $t\le T_m$, Lemmas~\ref{lem 4.14}--\ref{m19:moments} and
 \eqref{eq 6.1}  imply
\begin{align}\label{eq 6.12}
 |m[q\varpi]|&=\cO(\nu^{-1}\ve^\oo), \quad 
 |M[q\varpi]|\lesssim \nu \ve^2 \beta_\ve^{N_1} + \nu\ve\beta_\ve^N,\quad 
 |\dbz_p|\lesssim
\nu\ve\beta_\ve^{N_1+1}
                 +\nu\|q\|_{D_\ve}
                 +\nu^2\beta_\ve^N.
\end{align}
For $\nu$ small enough,  \eqref{eq 6.1} and Lemma~\ref{bd19:J1} give
\begin{equation}\label{eq 6.13}
    |\cJ_1|\le\ka\|q\|_{D_\ve}^2/8.
\end{equation}
By inserting \eqref{eq 6.12} into Proposition~\ref{prop 4.9}, we obtain
\begin{align}\label{eq 6.14}
 D_\ve[q]
 &\le-\ka\|q\|_{D_\ve}^2
       +C\bigl(\ve^2+\beta_\ve^{2N-2N_1}\bigr)
                         \nu^2\ve^2\beta_\ve^{2N_1}.
\end{align}
Moreover, by inserting  \eqref{eq 6.12} into Lemmas~\ref{lem 4.13}, \ref{lem 4.11}, \ref{lem 4.12},  we get
\begin{align}
\label{eq 6.15}
 |\cJ_2|&\le C\beta_\ve^{N-N_1}
                          \nu^2\ve^2\beta_\ve^{2N_1},
\\\label{eq 6.16}
 |\cJ_3|&\le C\ve\beta_\ve\|q\|_{D_\ve}^2
   +C\nu\ve^3\beta_\ve^{N+N_1}
       \left(\nu\ve\beta_\ve^{N_1+1}
                     +\nu\|q\|_{D_\ve}
                     +\nu^2\beta_\ve^N\right)\notag\\
 &\le\f\ka8\|q\|_{D_\ve}^2
   +C\left(\ve^2\beta_\ve^{N+1}
           +\nu\ve\beta_\ve^{2N-N_1}
           +\nu^2\ve^4\beta_\ve^{2N}\right)
                          \nu^2\ve^2\beta_\ve^{2N_1},
\\\label{eq 6.17}
 |\cJ_4|&\le C\ve\beta_\ve^{N_1}
                            \|q\|_{D_\ve}^2
      +C\ve^2\beta_\ve^{N_1}
       \left(\nu\ve\beta_\ve^{N_1+1}
                     +\nu\|q\|_{D_\ve}
                     +\nu^2\beta_\ve^N\right)
                                      \|q\|_{D_\ve}\notag\\
 &\le\f\ka8\|q\|_{D_\ve}^2
       +C\left(\ve^4\beta_\ve^{2N_1+2}
                        +\nu^2\ve^2\beta_\ve^{2N}\right)
                          \nu^2\ve^2\beta_\ve^{2N_1}. 
\end{align}

Combining \eqref{eq 6.13}–\eqref{eq 6.17} with \eqref{eq 4.24}, we obtain
\begin{equation}\label{eq 6.18}
 \br_a^2t\f{d}{dt}
          \left(\f{E_\ve[q]}{\br_a^2}\right)
       +\f\ka2\|q\|_{D_\ve}^2
 \lesssim \lam\nu^2\ve^2\beta_\ve^{2N_1},
\end{equation}
where one may take $N_1\gg1$ and $\nu\ll1$ so that
\begin{equation}\label{eq 6.19}
 \lam \eqdefa 
      \ve^2+\beta_\ve^{2N-2N_1} +\beta_\ve^{N-N_1} +\ve^2\beta_\ve^{N+1}
           +\nu\ve\beta_\ve^{2N-N_1}
           +\nu^2\ve^4\beta_\ve^{2N} +\ve^4\beta_\ve^{2N_1+2}
                        +\nu^2\ve^2\beta_\ve^{2N} \ll 1.
\end{equation}
Since $\br_a\approx r_0$ and $E_\ve[q](t)\rightarrow 0 $ as $t\rightarrow 0$, integrating \eqref{eq 6.18}  yields for any $t\leq T_m$
\begin{equation}\label{eq 6.20}
 E_\ve[q](t) + \int_0^{t} \|q\|_{D_\ve}^2 \frac{d\tau}{\tau}
 \lesssim \Big(\sup_{0<\tau\le T_m}\lam(\tau)\Big)
                          \nu^2\ve^2\beta_\ve^{2N_1}.
\end{equation}
By using Proposition~\ref{prop 4.2}, \eqref{eq 6.12} and the above inequality, we arrive at
\begin{align}\label{eq 6.21}
 \|q\|_{X_\ve}^2 + \int_0^{t} \|q\|_{D_\ve}^2 \frac{d\tau}{\tau} \lesssim   \Big(\sup_{0<\tau\le T_m}\lam(\tau)\Big)
                          \nu^2\ve^2\beta_\ve^{2N_1} 
       +\ve^2\|q\|_{X_\ve}^2
       +\nu^2\ve^2\beta_\ve^{2N} \le\f9{16}\nu^2\ve^2\beta_\ve^{2N_1}, 
\end{align}
which  implies $T_m=\nu^{-s}$ by a bootstrap argument.

\noindent
{\bf Step 3. Return to physical variables.}
By integrating the estimate for $\dot{\bar z}_p$ in \eqref{eq 6.12} and using \eqref{eq 6.21}, we
obtain
\begin{align}\label{eq 6.22}
 |\bz_p(t)|
 &\lesssim \nu t\ve\beta_\ve^{N_1+1}
       +\nu^2t\beta_\ve^N
       +\nu
          \left(\int_0^t\|q\|_{D_\ve}^2\f{\dd\tau}{\tau}\right)^{1/2}
          \left(\int_0^t\tau\dd\tau\right)^{1/2}\notag\\
 &\lesssim \ve^3\beta_\ve^{N_1+1}
              +\nu\ve^2\beta_\ve^N
              +\nu\ve^3\beta_\ve^{N_1},
\end{align}
from which we deduce that for some $N$ large enough,
\begin{equation}\label{eq 6.23}
 \f{|\br_a-r_0|+|\bz_p|}{\sqrt{\nu t}}
                    \le C\ve\beta_\ve^N.
\end{equation}
By using the above inequality, Proposition \ref{prop 3.1}, \eqref{eq 6.21} and H\"older inequality, we bound
\begin{equation}\label{eq 6.24}
 \int_{\HH}\left|
   \om_\tht(t,r,z)
   -\f{1}{4\pi\nu t}
      \exp\left(-\f{(r-r_0)^2+(z-\bz_a(t))^2}{4\nu t}\right)
                         \right|\dd r\dd z
 \le C_s\ve\beta_\ve^N,
 \qquad 0<t\le \nu^{-s}.
\end{equation}

Finally, \eqref{eq 3.1}, \eqref{eq 3.7} and property (a6) imply $\dot{\bar z}_a=\dot{\bar 
z}_0^{\EL}+O(\mu|\log\mu|^N)=\dot{\bar z}_0^{\EL}+O((\nu t)^{1/2-\eta})$ for every fixed 
$\eta>0$, and the integral definition gives $\bar z_a(0)=z_0$. This completes the proof.
\end{proof}

\appendix
\renewcommand{\theequation}{\thesection.\arabic{equation}}
\setcounter{equation}{0}

\section{Auxiliary proofs}\label{appA}
This Appendix contains the proofs of several auxiliary results used in the main text.

\subsection{Proof of Lemma \ref{lem 1.4}}\begin{proof}[Proof of Lemma \ref{lem 1.4}]
The coordinate definitions give
\[
t\dot\varepsilon=\varepsilon\left(\frac12-\frac{t\dot{\bar r}}{\bar r}\right),
\qquad t\partial_t\varpi=t\dot\varepsilon\xi_1,
\qquad \partial_1\Omega=\varepsilon\frac{\Omega}{\varpi}
\quad\text{on }\partial\HH_\varepsilon,
\]
where the last quotient is understood through its smooth extension to the axis. Since $\Omega$ vanishes on $\partial\HH_\varepsilon$, differentiation of its mass and first moments produces no moving-boundary contribution.
Integration by parts in the diffusion terms yields
\[
\begin{aligned}
\int_{\HH_\varepsilon}
\left(\cL\Omega+\varepsilon\partial_1(\Omega/\varpi)\right)\dd\xi
&=-2\int_\R\partial_1\Omega(t,-\varepsilon^{-1},\xi_2)\dd\xi_2,\\
\int_{\HH_\varepsilon}\xi_1
\left(\cL\Omega+\varepsilon\partial_1(\Omega/\varpi)\right)\dd\xi
&=-\frac12M_1[\Omega]-\varepsilon\int_{\HH_\varepsilon}\frac{\Omega}{\varpi}\dd\xi+\frac2\varepsilon\int_\R\partial_1\Omega(t,-\varepsilon^{-1},\xi_2)\dd\xi_2,\\
\int_{\HH_\varepsilon}\xi_2
\left(\cL\Omega+\varepsilon\partial_1(\Omega/\varpi)\right)\dd\xi
&=-\frac12M_2[\Omega]-2\int_\R\xi_2\partial_1\Omega(t,-\varepsilon^{-1},\xi_2)\dd\xi_2.
\end{aligned}
\]
The axis conditions on $\Phi$ also imply
\[
\begin{aligned}
\int_{\HH_\varepsilon}\{\Phi,\Omega/\varpi\}\dd\xi&=0,\\
\int_{\HH_\varepsilon}\xi_1\{\Phi,\Omega/\varpi\}\dd\xi
&=\int_{\HH_\varepsilon}\partial_2\Phi\,\Omega/\varpi\dd\xi,\\
\int_{\HH_\varepsilon}\xi_2\{\Phi,\Omega/\varpi\}\dd\xi
&=-\int_{\HH_\varepsilon}\partial_1\Phi\,\Omega/\varpi\dd\xi.
\end{aligned}
\]
Finally, $\{\varpi^2,\Omega/\varpi\}=2\varepsilon\partial_2\Omega$. Substituting these identities into \eqref{eq 2.10} and integrating the frame terms against $1$, $\xi_1$, and $\xi_2$ completes the proof of Lemma \ref{lem 1.4}.
\end{proof}

\subsection{Proof of Lemma \ref{lem 2.2}}
\begin{proof}[Proof of Lemma \ref{lem 2.2}]
The gradient estimates follow from the kernel argument in \cite{GS2024}. Set $s=\ve^2|\xi-\zeta|^2/(\varpi(\xi)\varpi(\zeta))$. The expansions of $F$ at zero and infinity yield
\begin{equation}\label{eq A.1}
   |sF'(s)|\le C\min(1,s^{-3/2}),\qquad
 |F(s)-2sF'(s)|\le C
 \begin{cases}1+|\log s|,&s\le1,\\s^{-3/2},&s\ge1.\end{cases} 
\end{equation}
 
Direct differentiation of the kernel yields
\begin{equation}\label{eq A.2}
    \begin{split}
      \varpi^{-1}(\xi)\p_{\xi_2}K
 &=\f{1}{\pi}\sqrt{\frac{\varpi(\zeta)}{\varpi(\xi)}}
       sF'(s)\f{\xi_2-\zeta_2}{|\xi-\zeta|^2},\\
 \varpi^{-1}(\xi) \p_{\xi_1}K
 &=\f{1}{\pi}\sqrt{\frac{\varpi(\zeta)}{\varpi(\xi)}}
       sF'(s)\f{\xi_1-\zeta_1}{|\xi-\zeta|^2}
   +\f{\ve}{4\pi\varpi(\xi)}
       \sqrt{\f{\varpi(\zeta)}{\varpi(\xi)}}
       \bigl(F(s)-2sF'(s)\bigr).  
    \end{split}
\end{equation}
Thus by inserting 
\[
 \sqrt{\f{\varpi(\zeta)}{\varpi(\xi)}}\le1+\sqrt{s}, \andf \f{\ve}{\varpi(\xi)}
       \sqrt{\f{\varpi(\zeta)}{\varpi(\xi)}}
 =\f{\sqrt{s}}{|\xi-\zeta|}
       \f{\varpi(\zeta)}{\varpi(\xi)},
\]
into \eqref{eq A.2}, and then considering $s\le1$ and $s\ge1$ separately, we obtain
\begin{equation}\label{eq A.3}
 \f{|\grad_\xi\phi(\xi)|}{\varpi(\xi)}
 \le C\int_{\HH_\ve}
             \f{|\eta(\zeta)|}{|\xi-\zeta|}\dd\zeta.
\end{equation}
The Hardy--Littlewood--Sobolev inequality, followed by weighted
H\"older's inequality, proves \eqref{eq 2.17}. Using
$\bxi\le\langle\zeta\rangle+|\xi-\zeta|$ and splitting
the integral in \eqref{eq A.3} at any radius $R>0$ gives
\[
 \|\bxi\varpi^{-1}\grad\phi\|_{L^\oo}
 \le C\|\eta\|_{L^1}+
 C\bigl(R^{1/2}\|\bxi\eta\|_{L^4}
       +R^{-1/2}\|\bxi\eta\|_{L^{4/3}}\bigr).
\]
Optimization in $R$, the Sobolev embedding $H^1\hookrightarrow L^4$
on a half-plane, and weighted H\"older's inequality prove
\eqref{eq 2.18}.

To prove \eqref{eq 2.19}, we use the bound $|F(s)|\lesssim1+\log_+(1/s)$, where
$\log_+a\eqdefa\max(0,\log a)$. Moreover,
\begin{equation}\label{eq A.4}
    \log_+ \frac{1}{s} \lesssim \beta_\ve +\log\bxi +\log\langle \zeta \rangle + \log_+\frac{1}{|\xi-\zeta|},
\end{equation}
which implies
\begin{equation}\label{eq A.5}
    \begin{split}
|\phi(\xi)| &\lesssim \sqrt{\varpi(\xi)} \int_{\HH_\ve} \Big(\beta_\ve + \Big|\log_+\frac{1}{|\xi-\zeta|}\Big| + \log \bxi  + \log \langle \zeta \rangle\Big) \sqrt{\varpi(\zeta)} |\eta(\zeta)| \dd\zeta\\
&\lesssim \sqrt{\varpi(\xi)}
 \bigl(\beta_\ve+\log\bxi\bigr)\|\bxi^2\eta\|_{L^2}.
    \end{split}
\end{equation}
This proves \eqref{eq 2.19}. We thus complete the proof of Lemma \ref{lem 2.2}.
\end{proof}

\subsection{Proof of Lemma \ref{lem 2.1}}
\begin{proof}[Proof of Lemma \ref{lem 2.1}]
We first prove the uniform mapping estimate $\phi=\BS_\ve[f]=\cO_{\cS_\star}(\beta_\ve)$ for $f=\cO_\cZ(1)$. Indeed, Lemma \ref{lem 2.2} gives the estimate for $\phi$ and $\grad\phi$.
Since $[\p_2,\BS_\ve]=0$, it suffices to estimate the $\p_1$-derivatives of $\phi$. By applying $\p_1^{n-2}$ to
\begin{equation}\label{eq A.6}
-\grad_\xi\cdot (\varpi^{-1}\grad_\xi \phi)=f   \implies  \p_1^2\phi
 =-\varpi f+\ve\varpi^{-1}\p_1\phi
                    -\p_2^2\phi,
\end{equation}
we get for $n\geq 3$,
\begin{equation}\label{eq A.7}
\begin{split}
     \p_1^n\phi &= -\p_1^{n-2}(\varpi f +\p_2^2 \phi)+\ve \p_1^{n-3} \p_1(\varpi^{-1}\p_1\phi)\\
     &= -\p_1^{n-2}(\varpi f+\p_2^2 \phi)- \ve \p_1^{n-3}(f+\varpi^{-1}\p_2^2\phi)\\
     &= -\p_1^{n-2}(\varpi f+\p_2^2 \phi) -\ve \p_1^{n-3}f -  \int_0^1 \lam^{n-3} (\p_1^{n-2}\p_2^2\phi) (-\ve^{-1}+\lam(\xi_1+\ve^{-1}),\xi_2) \dd\lam.
\end{split}
\end{equation}
Thus by using the first line for the case $\varpi\geq \f12$ and the last line for the case $\varpi \leq \f12$, wherein $$\langle -\ve^{-1}+\lam(\xi_1+\ve^{-1}) \rangle \lesssim \bxi,$$
we obtain the higher derivative estimates follow from an induction argument on the degree of  $\xi_1$-derivative.
Using this mapping property and the fact that $\bxi^a \lesssim \ve^J \bxi^{J+a}$ for $|\xi| \geq (16\ve)^{-1}$, we obtain for any $M\in\N$,
\begin{equation*}
    (1-\chi(2\ve\xi))\left(\BS_\ve-
                \sum_{m=0}^N\ve^m\BS_m\right)[f],\quad  \left(\BS_\ve-\sum_{m=0}^N\ve^m\BS_m\right)[(1-\chi(2\ve\xi))f] = \cO_{\cS_\star}(\ve^M).
\end{equation*}
Here  $\chi(\zeta)$ is a smooth cut-off function satisfying $\chi=1$ for $|\zeta|\leq \f12$ and  $\chi=0$ for $|\zeta|\geq \f34$. 

It therefore remains to prove the following estimate in the region $|\ve\xi|,|\ve\zeta|\leq3/8$:
    \begin{equation}\label{eq A.8}
\chi(2\ve \xi) \left(\BS_\ve-
                \sum_{m=0}^N\ve^m\BS_m\right)[\chi(2\ve\xi)f]
  =\cO_{\cS_\star}
                 \bigl(\beta_\ve\ve^{N+1}\bigr).
\end{equation}
For $|\ve\xi|,|\ve\zeta|\leq 3/8$ and $\xi\ne\zeta$, we have
\[
s=\frac{\ve^2|\xi-\zeta|^2}{\varpi(\xi)\varpi(\zeta)}
\leq \frac{36}{25}<4,
\qquad
\log\frac{8}{\sqrt{s}}
=\beta_\ve+L+\frac12\log\bigl(\varpi(\xi)\varpi(\zeta)\bigr).
\]
The convergent expansion of $F$ therefore expresses $K$ as
$(\beta_\ve+L)$ times an analytic coefficient plus another
analytic coefficient, both analytic in $(\ve\xi,\ve\zeta)$
on a neighbourhood of this region, including the diagonal.
Taylor's formula with integral remainder for these coefficients,
together with
$(\p_{\xi_j}+\p_{\zeta_j})L=0$ for $j=1,2$,
gives, for some $p\in\N$,
\begin{equation}\label{eq A.9}
\begin{aligned}
 &\left|(\p_{\xi_1}+\p_{\zeta_1})^{\ga_1}(\p_{\xi_2}+\p_{\zeta_2})^{\ga_2}
 \left(K(\xi,\zeta)-\sum_{m=0}^N\ve^m
                   \bigl((\beta_\ve+L)P_m+Q_m\bigr)\right)\right|\\
 &\hspace{18mm}\leq C_{N,\ga}\ve^{N+1}
       (1+\beta_\ve+|L|)\bxi^p\langle\zeta\rangle^p.
\end{aligned}
\end{equation}
We then get by an induction argument (where we may need to modify the support of $\chi$ slightly in each step), Hilbert–Schmidt Lemma and integration by parts that for some $p'>p$,
\begin{equation}\label{eq A.10}
    \chi(2\ve \xi) \left|\grad_\xi^\ga
 \left(\BS_\ve-\sum_{m=0}^N\ve^m\BS_m\right)[\chi(2\ve\xi)f]\right|\lesssim_{N,\ga}\ve^{N+1} \beta_\ve \bxi^{p'} \sum_{|\al|\leq |\ga|}\|\bxi^{p'}\grad_\xi^\al f \|_{L^2}.
\end{equation}
We thus prove \eqref{eq 2.20} and in particular the estimate of the first term in the left-hand-side of \eqref{eq 2.21} when $\ga=0$. Notice that for $|\xi|,|\zeta|\leq 2 \ve^{-\sig_1}$, 
\begin{equation*}
    \Big|\varpi^{-1}\grad K -\grad \big( (\beta_\ve+L)P_0+Q_0\big)\Big| \lesssim \ve (\beta_\ve +|L|),
\end{equation*}
which together with Hilbert-Schmidt Lemma implies the estimate of second term in the left-hand-side of \eqref{eq 2.21}. We thus complete the proof of Lemma \ref{lem 2.1}.
\end{proof}

\subsection{Proof of Proposition \ref{en19:diffusion-prop}}
\begin{proof}[Proof of Proposition \ref{en19:diffusion-prop}]
We first split $D_\ve[q]=\cD_{\ve,1}[q]+\cD_{\ve,2}[q]+\cD_{\ve,3}[q]$,
where
\begin{equation}\label{eq A.11}
\begin{split}
 \cD_{\ve,1}[q]&\eqdefa
  \langle(\cL+3\ve\varpi^{-1}\p_1)q,qW_\ve\rangle
           +\f12\langle q^2,t\p_tW_\ve\rangle,\\
 \cD_{\ve,2}[q]&\eqdefa
  \langle(\cL+3\ve\varpi^{-1}\p_1)q,-\cB_\ve[q]\rangle,\\
  \cD_{\ve,3}[q]&\eqdefa -\f12\langle q,[t\pt,\cB_\ve]q\rangle.
\end{split}
\end{equation}
We first  show that 
\begin{align}
 \cD_{\ve,1}[q]
 =&-\int_{\HH_\ve}W_\ve|\grad_\xi q|^2\dd\xi
      -\int_{\HH_\ve}
        (\grad_\xi\log H_\ve\cdot\grad_\xi q)qW_\ve\dd\xi \notag\\
 &-\int_{\HH_\ve}
    \left(W_\ve V_\ve+\varpi^3
       +\f{t\dbr_a}{2\br_a}\ve\xi_1H_\ve\right)q^2\dd\xi-\ve\int_{\HH_\ve}(\p_1H_\ve)q^2\dd\xi
  -\ve\int_{\p\HH_\ve}H_\ve q^2\dd\xi_2, \label{eq A.12}\\
  \cD_{\ve,2}[q] =&\int_{\HH_\ve}\varpi^3q^2\dd\xi
         +\ve\zeta_\ve[q], \label{eq A.13} \\
    \cD_{\ve,3}[q]=&-2t\pt\ve\,\zeta_\ve[q]  \label{eq A.14}.
\end{align}

\noindent {\bf Step 1. Proof of \eqref{eq A.12}.}
Since $W_\ve=\varpi H_\ve$ and $W_\ve=0$ on
$\p\HH_\ve$, by integration by parts, we get
\begin{align*}
 &\langle\lap_\xi q,qW_\ve\rangle
  =-\int_{\HH_\ve}W_\ve|\grad_\xi q|^2\dd\xi
    -\int_{\HH_\ve}(\grad_\xi W_\ve\cdot\grad_\xi q)q\dd\xi,\\
& \left\langle\f12\xi\cdot\grad_\xi q+q,qW_\ve\right\rangle
  =-\int_{\HH_\ve}
      \left(\f14\xi\cdot\grad_\xi W_\ve
                          -\f12W_\ve\right)q^2\dd\xi,\\
&\int_{\HH_\ve}H_\ve q\p_1q\dd\xi
 =-\f12\int_{\HH_\ve}(\p_1H_\ve)q^2\dd\xi
  -\f12\int_{\p\HH_\ve}H_\ve q^2\dd\xi_2.
\end{align*}
And then by using 
\begin{align*}
 &\grad_\xi W_\ve
  =W_\ve\grad_\xi\log H_\ve
                              +\ve H_\ve(1,0)^{\rm T},\\
 &\f14\xi\cdot\grad_\xi W_\ve-\f12W_\ve
                       -\f12t\p_tW_\ve
  =W_\ve V_\ve+\varpi^3
      +\left(\f\ve4-\f12t\pt\ve\right)\xi_1H_\ve=W_\ve V_\ve+\varpi^3
      +\f{t\dbr_a}{2\br_a}\ve\xi_1H_\ve,
\end{align*}
we complete the proof of \eqref{eq A.12}.

\noindent {\bf Step 2. Proof of \eqref{eq A.13}.}
Now for \eqref{eq A.13}, recall that
\begin{equation}\label{eq A.15}
 -\grad_\xi\cdot(\varpi^{-1}\grad_\xi\phi)=\varpi q,
\end{equation}
which gives 
\begin{align*}
 &-\int_{\HH_\ve}\varpi\phi
                      (\lap_\xi q+3\ve\varpi^{-1}\p_1q)\dd\xi=\int_{\HH_\ve}\varpi\grad_\xi\phi\cdot\grad_\xi q\dd\xi
               -2\ve\int_{\HH_\ve}\phi\p_1q\dd\xi\\
 &\quad=\int_{\HH_\ve}\varpi q
                       (-\lap_\xi\phi+\ve\varpi^{-1}\p_1\phi)\dd\xi
       =\int_{\HH_\ve}\varpi^3q^2\dd\xi.
\end{align*}
For the remaining part of $\cL$, integration by parts yields
\begin{equation}\label{eq A.16}
\begin{split}
&-\int_{\HH_\ve}\varpi\phi
                  \left(\f12\xi\cdot\grad_\xi q+q\right)\dd\xi
 =\f12\int_{\HH_\ve}\varpi q\,\xi\cdot\grad_\xi\phi\dd\xi
   +\f\ve2\int_{\HH_\ve}\xi_1\phi q\dd\xi\\
 &=-\f12\int_{\HH_\ve}
    \grad_\xi\cdot(\varpi^{-1}\grad_\xi\phi)\,\xi\cdot\grad_\xi\phi\dd\xi  +\f\ve2\int_{\HH_\ve}\xi_1\phi q\dd\xi\\
 &=\f12\int_{\HH_\ve}\f{|\grad_\xi\phi|^2}{\varpi}\dd\xi
   +\f14\int_{\HH_\ve}
          \f{\xi\cdot\grad_\xi|\grad_\xi\phi|^2}{\varpi}\dd\xi +\f\ve2\int_{\HH_\ve}\xi_1\phi q\dd\xi\\
 &=\f\ve4\int_{\HH_\ve}
             \f{\xi_1|\grad_\xi\phi|^2}{\varpi^2}\dd\xi +\f\ve2\int_{\HH_\ve}\xi_1\phi q\dd\xi= \ve \zeta_\ve[q].
\end{split}
\end{equation}

\noindent {\bf Step 3. Proof of \eqref{eq A.14}.}
Finally for the calculation of $\cD_{\ve,3}$, since $\cB_\ve[q](\xi)=\varpi(\xi)\phi(\xi)$, we write
\begin{equation}\label{eq A.17}
 E_{\cB_\ve}(\ve)\eqdefa\langle q,\cB_\ve[q]\rangle
 =\iint_{\HH_\ve\times\HH_\ve}
              B(\xi,\zeta)q(\xi)q(\zeta)\dd\xi\dd\zeta,
\end{equation}
where 
\begin{equation*}
 B(\xi,\zeta)\eqdefa  
 \f1{2\pi}\bigl(\varpi(\xi)\varpi(\zeta)\bigr)^{3/2}
 F\Bigl(\frac{\ve^2|\xi-\zeta|^2}{\varpi(\xi)\varpi(\zeta)}\Bigr),
\end{equation*}
from which we infer
\begin{equation}\label{eq A.18}
    \begin{split}
\ve\p_\ve B(\xi,\zeta) = (\xi\cdot\grad_\xi + \zeta\cdot\grad_\zeta) B(\xi,\zeta).
    \end{split}
\end{equation}
Thus on the one hand, 
\begin{equation}\label{eq A.19}
 \p_\ve E_{\cB_\ve}(\ve)
 =\iint_{\HH_\ve\times\HH_\ve}
       \p_\ve B(\xi,\zeta)q(\xi)q(\zeta)\dd\xi\dd\zeta.
\end{equation}
While on the other hand, by \eqref{eq A.16} and symmetrization, we can verify that 
\begin{equation}\label{eq A.20}
    \begin{split}
4\zeta_\ve[q] &= \int_{\HH_\ve}
\frac{\xi_1|\nabla_\xi\phi|^2}
{\varpi(\xi)^2}\,d\xi
+2\int_{\HH_\ve}\xi_1\phi q\,d\xi = \frac{2}{\ve}\int_{\HH_\ve} \eta \xi\cdot \grad \phi d\xi + 2 \int_{\HH_\ve}\xi_1\phi q\,d\xi \\
&= \frac{1}{\ve} \iint_{\HH_\ve\times\HH_\ve}  (\xi\cdot\grad_\xi + \zeta\cdot\grad_\zeta)\Big(\frac{1}{\varpi(\xi)\varpi(\zeta)}B(\xi,\zeta)\Big) \varpi(\xi)\varpi(\zeta) q(\xi)q(\zeta) \dd\xi\dd\zeta   \\
&\quad +\iint_{\HH_\ve\times\HH_\ve}  \Big(\frac{\xi_1}{\varpi(\xi)}+\frac{\zeta_1}{\varpi(\zeta)}\Big)\Big(\frac{1}{\varpi(\xi)\varpi(\zeta)}B(\xi,\zeta)\Big) \varpi(\xi)\varpi(\zeta) q(\xi)q(\zeta) \dd\xi\dd\zeta\\
& =\frac{1}{\ve}\iint_{\HH_\ve\times\HH_\ve}
(\xi\cdot\grad_\xi + \zeta\cdot\grad_\zeta) B(\xi,\zeta)
q(\xi)q(\zeta)\,d\xi\,d\zeta.
    \end{split}
\end{equation}

Then in view of \eqref{eq A.18}-\eqref{eq A.20}, we obtain
\begin{equation}\label{eq A.21}
    \begin{split}
\langle q, [t\pt, \cB_\ve][q] \rangle = \langle q, t\pt \cB_\ve  [q] \rangle = t\pt\ve \times  \partial_\ve E_{\mathcal B_\ve}(\ve) = 4 t\pt\ve  \zeta_\ve[q],
    \end{split}
\end{equation}
which proves \eqref{eq A.14}. 

\noindent {\bf Step 4. Proof of \eqref{eq 4.33}.} Finally by rearranging the decomposition of $D_\ve[q]$ and then using $t\p_t \ve =\ve(\f12-\frac{t\dbr_a}{\br_a})$, we complete the proof of Proposition \ref{en19:diffusion-prop}.
\end{proof}

\section*{Acknowledgements}

P. Zhang is partially supported by the National Key R\&D Program of China under grant 2021YFA1000800 and by the National Natural Science Foundation of China under grants 12421001, 12494542, and 12288201. The authors used OpenAI models as assistive tools in
preparing this manuscript. Mathematical arguments, proofs, and computations in the manuscript
were verified and finalized by the authors.

\end{document}